\documentclass[12pt]{article}                
\usepackage{graphicx}
\usepackage{psfrag}
\usepackage{amsmath, amssymb}
\usepackage{mathptmx}
\usepackage{ntheorem}

\newcommand{\A}{{\cal A}}
\newcommand{\T}{{\cal T}}
\newcommand{\M}{{\cal M}}

\newcommand{\DD}{{\cal D}}
\newcommand{\PP}{{\cal P}}
\newcommand{\EE}{{\cal E}}
\newcommand{\RR}{{\cal R}}
\newcommand{\LL}{{\cal L}}

\newcommand{\CC}{{\cal C}}
\newcommand{\OO}{{\cal O}}
\newcommand{\OC}{{{\cal O}({\cal C})}}
\newcommand{\WW}{{\cal W}}
\newcommand{\WC}{{{\cal W}({\cal C})}}
\newcommand{\Ind}{{\hspace{0.3mm}{\rm I}\hspace{0.1mm}}}
\newcommand{\eps}{{\varepsilon}}

\newcommand{\smsp}{\hspace{0.3mm}}
\newcommand{\e}{\mathbb{E}}

\newcommand{\p}{\mathbb{P}}

\newcommand{\Reals}{\mathbb{R}}
\newcommand{\Natural}{\mathbb{N}}
\newcommand{\la}{\langle}
\newcommand{\ra}{\rangle}
\newcommand\qed{\hfill\hbox{\rlap{$\sqcap$}$\sqcup$}}
\newcommand{\bs}{{\bar{\sigma}}}

\newtheorem{lemma}{Lemma}
\newtheorem{theorem}{Theorem}
\newtheorem{corollary}{Corollary}

\theoremstyle{nonumberplain}
\newcommand\specialref{}
\newtheorem{thmx}{Theorem \specialref}
\newenvironment{thmref}[2][$'$]
  {\renewcommand\specialref{\ref{#2}#1}\thmx}
  {\endthmx}

\begin{document}

\title{Hierarchical exchangeability of pure states\\ in mean field spin glass models}
\author{Dmitry Panchenko\thanks{Dept. of Mathematics, Texas A\&M University, panchenk@math.tamu.edu. Partially supported by NSF grant.}\\
}
\date{}
\maketitle
\begin{abstract}
The main result in this paper is motivated by the M\'ezard-Parisi ansatz which predicts a very special structure for the distribution of spins in diluted mean field spin glass models, such as the random $K$-sat model at positive temperature. Using the fact that one can safely assume the validity of the Ghirlanda-Guerra identities in these models, we prove hierarchical exchangeability of pure states for the asymptotic Gibbs measures, which allows us to apply a representation result for hierarchically exchangeable arrays recently proved in \cite{AP}. Comparing this representation with the predictions of the M\'ezard-Parisi ansatz, one can see that the key property still missing is that the multi-overlaps between pure states depend only on their overlaps.

\end{abstract} 
\vspace{0.5cm}
Key words: spin glasses, diluted models, exchangeability.\\
Mathematics Subject Classification (2010): 60K35, 60G09, 82B44

\section{Introduction}

Many mean field spin glass models are described by a random Hamiltonian $H_N(\sigma)$ on the space of spin configurations $\Sigma_N = \{-1,+1\}^N$ (see \cite{MPV}, \cite{SG} or \cite{SG2}). For example, in the classical Sherrington-Kirkpatrick model \cite{SK}, 
\begin{equation}
H_N(\sigma) = \frac{1}{\sqrt{N}}\sum_{i,j=1}^N g_{i,j} \sigma_i \sigma_j,
\label{Ham1}
\end{equation}
where $(g_{i,j})_{i,j\geq 1}$ are i.i.d. standard Gaussian random variables, while in the random $K$-sat model,
\begin{equation}
H_N(\sigma) = \sum_{k\leq \pi(\alpha N)} \prod_{1\leq j\leq K} \frac{1+\eps_{j,k} \sigma_{i_{j,k}}}{2},
\label{Ham2}
\end{equation}
where $\alpha>0$ is called the connectivity parameter, $\pi(\alpha N)$ is a Poisson random variable with the mean $\alpha N$, $(\eps_{j,k})_{j,k\geq 1}$ are independent Rademacher random variables, and the indices $(i_{j,k})_{j,k\geq 1}$ are independent uniform on $\{1,\ldots,N\}$. The random $K$-sat model is an example of a so called diluted model (diluted refers to the fact that the (hyper-)graph of interactions between spins is sparse), and the main goal of this paper is to make some progress toward the M\'ezard-Parisi ansatz for diluted models originating in \cite{Mezard}. The reason the above two models are called mean field models is because the distributions of their Hamiltonians are invariant under the permutations of coordinates $\sigma_1,\ldots,\sigma_N$. This property is called \emph{symmetry between sites}.

The main goal in spin glass models is usually to compute the limit of the free energy 
\begin{equation}
F_N = \frac{1}{N}\e \log \sum_{\sigma\in \Sigma_N} \exp \bigl(-\beta H_N(\sigma) \bigr)
\end{equation}
as $N\to\infty$, for all inverse temperature parameters $\beta>0$. In the Sherrington-Kirkpatrick model, the formula for the free energy was famously invented by Parisi in \cite{Parisi79,Parisi} and proved rigorously by Talagrand in \cite{TPF} following important work of Guerra in \cite{Guerra}, who showed that the Parisi formula is an upper bound on the free energy. A more recent proof of the Parisi formula in \cite{PPF} was based on understanding the structure of the Gibbs measure in the infinite-volume limit predicted by the physicists in the eighties (see \cite{MPV}; this direction of research was jump-started in \cite{AA}). 

For diluted models, like the random $K$-sat model, the analogue of the Parisi formula for the free energy was proposed by M\'ezard and Parisi in \cite{Mezard} in the so called $1$-RSB case (a replica symmetric solution was proposed earlier in \cite{MZ}), but it was also clear what the natural extension of their solution should look like in the general case. A detailed description of the general M\'ezard-Parisi formula can be found, for example, in \cite{PT}. The fact that this formula gives an upper bound on the free energy was proved by Franz and Leone in \cite{FL}, which was the analogue of Guerra's bound \cite{Guerra} in the SK model. One approach to proving the matching lower bound was given in \cite{Pspins} (see Theorem $2$ and, in particular, Section 2.2 there), where the problem was reduced via an analogue of the Aizenman-Sims-Starr scheme \cite{AS2} to showing that the structure of the Gibbs measure in the infinite-volume limit is described by the functional order parameter proposed by M\'ezard and Parisi in \cite{Mezard}. Our main result here combined with the main result in \cite{AP} can be viewed as a partial progress in this direction, and at the end of the introduction we will explain what the remaining gap is. Let us mention that the M\'ezard-Parisi ansatz has been proved in full generality in the setting of the Sherringon-Kirkpatrick model and $p$-spin models (see Chapter 4 in \cite{SKmodel}), but the proof heavily relies on the special Gaussian nature of the Hamiltonian (\ref{Ham1}). In diluted models, where this ansatz is of real interest, the problem is still open in general, with one special case handled recently in \cite{1RSB}. 

In this paper, we will not work with any particular model and will simply assume that the asymptotic Gibbs measures satisfy the Ghirlanda-Guerra identities \cite{GG}. In the next section we will review how the Ghirlanda-Guerra identities arise in spin glass models and, as an example, show that one can safely assume their validity in the random $K$-sat model. The Ghirlanda-Guerra identities will be stated in this paper in a slightly more general form than usual to accommodate the more general notion of the asymptotic Gibbs measures in models other than the SK model but, of course, one gets this more general form for free from the usual proof of these identities.

Let us begin by recalling the definition of asymptotic Gibbs measures introduced in \cite{Pspins} (see also \cite{Austin2} for a different approach via exchangeable random measures). The Gibbs measure $G_N$ corresponding to the Hamiltonian $H_N(\sigma)$ is a (random) probability measure on $\{-1,+1\}^N$ defined by 
\begin{equation}
G_N(\sigma) = \frac{1}{Z_N} \exp\bigl(- \beta H_N(\sigma)\bigr)
\label{GibbsN}
\end{equation} 
where the normalizing factor $Z_N$ is called the partition function. Let $(\sigma^\ell)_{\ell\geq 1}$ be an i.i.d. sequence of replicas from the Gibbs measure $G_N$ and let $\mu_N$ be the joint distribution of the array of all spins on all replicas $(\sigma_i^\ell)_{1\leq i\leq N, \ell\geq 1}$ under the average product Gibbs measure $\e G_N^{\otimes \infty}$,
\begin{equation}
\mu_N\Bigl( \bigl\{\sigma_i^\ell = a_i^\ell \ :\ 1\leq i\leq N, 1\leq \ell \leq n \bigr\} \Bigr)
=
\e G_N^{\otimes n}\Bigl( \bigl\{\sigma_i^\ell = a_i^\ell \ :\ 1\leq i\leq N, 1\leq \ell \leq n \bigr\} \Bigr)
\label{muN}
\end{equation}
for any $n\geq 1$ and any $a_i^\ell \in\{-1,+1\}$.  We extend $\mu_N$ to a distribution on $\{-1,+1\}^{\Natural\times\Natural}$ by setting $\sigma_i^\ell=1$ for $i\geq N+1.$ Let $\M$ be the sets of all possible limits of $(\mu_N)$ over subsequences with respect to the weak convergence of measures on the compact product space $\{-1,+1\}^{\Natural\times\Natural}$. Because of the symmetry between sites in mean field models, these measures inherit from $\mu_N$ the invariance under the permutation of both spin and replica indices $i$ and $\ell.$ By the Aldous-Hoover representation \cite{Aldous, Hoover2}, for any $\mu\in \M$, there exists a measurable function $s:[0,1]^4\to\{-1,+1\}$ such that $\mu$ is the distribution of the array 
\begin{equation}
s_i^\ell=s(w,u_\ell,v_i,x_{i,\ell}),
\label{sigma}
\end{equation}
where the random variables $w,(u_\ell), (v_i), (x_{i,\ell})$ are i.i.d. uniform on $[0,1]$. The function $s$ is defined uniquely for a given $\mu\in \M$ up to measure-preserving transformations (Theorem 2.1 in \cite{Kallenberg}), so we can identify the distribution $\mu$ of array $(s_i^\ell)$ with $s$. Since $s$ takes values in $\{-1,+1\}$, the distribution $\mu$ can actually be encoded by the function
\begin{equation}
{\sigma}(w,u,v) = \e_x\, s(w,u,v,x)
\label{fop}
\end{equation}
where $\e_x$ is the expectation in $x$ only. The last coordinate $x_{i,\ell}$ in (\ref{sigma}) is independent for all pairs $(i,\ell)$, so it plays the role of ``flipping a coin" with the expected value $\sigma(w,u_\ell,v_i)$. In fact, given the function (\ref{fop}), we can, obviously, redefine $s$ by
\begin{equation}
s(w,u_\ell,v_i,x_{i,\ell}) = 2 \Ind \Bigl(x_{i,\ell} \leq \frac{1+ \sigma(w,u_\ell,v_i) }{2}\Bigr) -1
\label{sigmatos}
\end{equation}
without affecting the distribution of the array $(s_i^\ell)$. This allows us to separate the randomness of the last coordinate $x_{i,\ell}$ from the randomness of the array $(\sigma(w,u_\ell,v_i))$ generated by the function $\sigma(w,u,v)$. 

Then we change the perspective as follows. Let $du$ and $dv$ denote the Lebesgue measure on $[0,1]$ and let us define a (random) probability measure 
\begin{equation}
G = G_w = du \circ \bigl(u\to \sigma(w,u,\cdot)\bigr)^{-1}
\label{Gibbsw}
\end{equation}
on the space of functions of $v\in [0,1]$,
\begin{equation}
H = \bigl\{ \|\sigma\|_\infty \leq 1 \bigr\}
\end{equation} 
(the unit ball of $L^\infty$), equipped with the topology of $L^2([0,1], dv)$. We will denote by $\sigma^1\cdot \sigma^2$ the scalar product in $L^2([0,1], dv)$ and by $\|\sigma\|$ the corresponding $L^2$ norm. The random measure $G$ in (\ref{Gibbsw}) is what we call the \emph{asymptotic Gibbs measure}, which encodes the limit $\mu\in \M$ above. The whole process of generating spins according to $\mu\in \M$ can now be visualized in several steps. First, we generate the Gibbs measure $G=G_w$ using the uniform random variable $w$. An i.i.d. sequence $\sigma^\ell = \sigma(w,u_{\ell},\cdot)$ for $\ell\geq 1$  of replicas from $G$ gives us a sequence of functions in $H$. Then, we plug in i.i.d. uniform random variables $(v_i)_{i\geq 1}$ into these functions to obtain the array $\sigma^\ell(v_i) = \sigma(w,u_\ell,v_i)$ and, finally, use it to generate spins as in (\ref{sigmatos}). From now on, we will keep the dependence of $G$ on $w$ implicit, denote i.i.d. replicas from $G$ by $(\sigma^\ell)_{\ell\geq 1}$ (which are now functions on $[0,1]$) and no longer explicitly use the random variables $(u_{\ell})$, and denote the sequence of spins (\ref{sigmatos}) corresponding to the replica $\sigma^\ell$ by
\begin{equation}
S(\sigma^\ell) = \Bigl( 2 \Ind\Bigl(x_{i,\ell} \leq \frac{1+ \sigma^\ell(v_i) }{2}\Bigr) -1 \Bigr)_{i\geq 1} \in \{-1,+1\}^\Natural.
\label{SpinsEll}
\end{equation}
Given $n\geq 1$ and replicas $\sigma^1,\ldots, \sigma^n$, we will denote the array of spins corresponding to these replicas by 
\begin{equation}
S^n = \bigl(S(\sigma^\ell) \bigr)_{1\leq \ell\leq n}.
\label{Sn}
\end{equation}
We will denote by $\la\cdot\ra$ the average with respect to $G^{\otimes \infty}$ (corresponding to the average in $(u_\ell)_{\ell\geq 1}$ in the sequence $(\sigma(w,u_{\ell},\cdot))_{\ell\geq 1}$) and by $\e$ the expectation with respect to $w$, $(v_i)$ and $(x_{i,\ell})$. In the definition of $\la\cdot\ra$ one can also include averaging in the random variables $(x_{i,\ell})$, since they depend on the replica index $\ell$, and such convention would be especially necessary if we dealt with the cavity computations (see e.g. \cite{Pspins, 1RSB}), when averaging in spins $S(\sigma^\ell)$ can also appear in the denominator. However, throughout this paper this will not happen and, by the linearity of expectation, we can think of averaging in $(x_{i,\ell})$ as a part of the expectation $\e$.

Because of the geometric nature of the asymptotic Gibbs measures $G$ as measures on the subset of $L^2([0,1],dv)$, the distance and scalar product between replicas play a crucial role in the description of the structure of $G$. We will denote the scalar product between replicas $\sigma^\ell$ and $\sigma^{\ell'}$ by $R_{\ell,\ell'} = \sigma^\ell\cdot \sigma^{\ell'}$, which is more commonly called \emph{the overlap} of $\sigma^\ell$ and $\sigma^{\ell'}$. For $\ell\ne\ell'$, let us notice that the overlap $R_{\ell,\ell'}$ is a function of spin sequence (\ref{SpinsEll}) generated by $\sigma^\ell$ and $\sigma^{\ell'}$ since, by the strong law of large numbers,
\begin{equation}
R_{\ell,\ell'} = \int \! \sigma^\ell(v) \sigma^{\ell'}(v)\, dv =
\lim_{j\to\infty} \frac{1}{j}\sum_{i=1}^j S\bigl(\sigma^\ell\bigr)_i \,S\bigl(\sigma^{\ell'}\bigr)_i
\label{overlapinfty}
\end{equation}
almost surely. We mention this here just to emphasize an obvious point that the array $S^n$ in (\ref{Sn}) contains much more information about the replicas on the space $H$ than just their overlaps. For example, one can similarly compute the multi-overlaps between replicas with pairwise distinct indices $\ell_1,\ldots,\ell_n$,
\begin{equation}
R_{\ell_1,\ldots,\ell_n} = \int \! \sigma^{\ell_1}(v) \cdots \sigma^{\ell_n}(v)\, dv =
\lim_{j\to\infty} \frac{1}{j}\sum_{i=1}^j S\bigl(\sigma^{\ell_1}\bigr)_i \cdots S\bigl(\sigma^{\ell_n}\bigr)_i.
\label{multioverlapinfty}
\end{equation}

From now on we will assume that the measure $G$ satisfies \emph{the Ghirlanda-Guerra identities}, which means that for any $n\geq 2,$ any bounded measurable function $f$ of the spins $S^n$ in (\ref{Sn}) and any bounded measurable function $\psi$ of one overlap,
\begin{equation}
\e \bigl\la f(S^n)\psi(R_{1,n+1}) \bigr\ra = \frac{1}{n}\hspace{0.3mm} \e\bigl\la f(S^n) \bigr\ra \hspace{0.3mm} \e\bigr\la \psi(R_{1,2})\bigr\ra + \frac{1}{n}\sum_{\ell=2}^{n}\e\bigl\la f(S^n)\psi(R_{1,\ell})\bigr\ra.
\label{GG}
\end{equation}
Another way to express the Ghirlanda-Guerra identities is to say that, conditionally on $S^n$, the law of $R_{1,n+1}$ is given by  the mixture 
\begin{equation}
\frac{1}{n} \hspace{0.3mm}\zeta + \frac{1}{n}\hspace{0.3mm} \sum_{\ell=2}^n \delta_{R_{1,\ell}},
\label{GGgen}
\end{equation}
where $\zeta$ denotes the distribution of $R_{1,2}$ under the measure $\e G^{\otimes 2}$,
\begin{equation}
\zeta(\ \cdot\ ) = \e G^{\otimes 2}\bigl(R_{1,2}\in\ \cdot\ \bigr).
\label{zeta}
\end{equation}
The identities (\ref{GG}) are usually proved for the function $f$ of the overlaps $(R_{\ell,\ell'})_{\ell,\ell'\leq n}$ instead of $S^n$, but exactly the same proof yields (\ref{GG}) as well (see e.g. Section 3.2 in \cite{SKmodel}). It is well known that these identities arise from the Gaussian integration by parts of a certain Gaussian perturbation Hamiltonian against the test function $f$, and one is free to choose this function to depend on all spins and not only overlaps.

In this paper we will be interested to say something about the distribution of the array of spins generated by the Gibbs measure $G$, but if one is only interested in the behavior of the overlaps then it is now known that the Ghirlanda-Guerra identities completely describe the measure in this sense in terms of the \emph{functional order parameter} $\zeta$ in (\ref{zeta}). Let us first list several purely geometric consequences.
\begin{enumerate}

\item[(i)] (\cite{SG} or Theorem 2.16 in \cite{SKmodel}) By Talagrand's positivity principle, the overlaps can take only nonnegative values, $\zeta([0,\infty))=1$.

\item[(ii)] (\cite{PGG} or Theorem 2.15 in \cite{SKmodel}) With probability one over the choice of random measure $G$ the following holds. If $q^*$ is the largest point in the support $\mbox{supp}(\zeta)$ of measure $\zeta$ then $G(\sigma : \|\sigma\|^2 = q^*)=1$. If $\zeta(\{q^*\})>0$ then $G$ is purely atomic, otherwise, $G$ has no atoms. 

\item[(iii)] (\cite{PUltra} or Theorem 2.14 in \cite{SKmodel}) With probability one, the support of $G$ is ultrametric, i.e. $G^{\otimes 3}(R_{2,3} \geq \min(R_{1,2},R_{1,3}))=1$.  

\end{enumerate}
When $G$ is purely atomic, its atoms are called \emph{pure states}. Otherwise, we will define pure states in some approximate sense that will be explained more precisely below. By ultrametricity, for any $q\in[0,q^*]$, the relation defined by
\begin{equation}
\sigma\sim_q \sigma' \Longleftrightarrow \sigma\cdot\sigma' \geq q
\label{qclusters}
\end{equation}
is an equivalence relation on the support of $G$. We will call these $\sim_q$ equivalence clusters simply $q$-clusters. Throughout the paper we will use the convention that, whenever we write $\sigma$, it belongs to the support of $G$ rather than the ambient space $H$.

To state our main result, let us first describe what is called the \emph{$r$-step replica symmetry breaking} (RSB) approximation, which means that we will group the values of the overlap into $r+1$ groups. Let us consider integer $r\geq 1$ that will be fixed throughout the paper. Consider an infinitary rooted tree of depth $r$ with the vertex set
\begin{equation}
\A = \Natural^0 \cup \Natural \cup \Natural^2 \cup \ldots \cup \Natural^r,
\label{Atree}
\end{equation}
where $\Natural^0 = \{*\}$, $*$ is the root of the tree and each vertex $\alpha=(n_1,\ldots,n_p)\in \Natural^{p}$ for $p\leq r-1$ has children 
$$
\alpha n : = (n_1,\ldots,n_p,n) \in \Natural^{p+1}
$$
for all $n\in \Natural$. Each vertex $\alpha$ is connected to the root $*$ by the path
$$
* \to n_1 \to (n_1,n_2) \to\cdots\to (n_1,\ldots,n_p) = \alpha.
$$
We will denote the set of vertices in this path (excluding the root) by 
\begin{equation}
p(\alpha) = \bigl\{n_1, (n_1,n_2),\ldots,(n_1,\ldots,n_p)  \bigr\}.
\label{pathtoleaf}
\end{equation}
We will denote by $|\alpha|$ the distance of $\alpha$ from the root (the same as cardinality of $p(\alpha)$). We will write $\alpha \succ \beta$ if $\beta \in p(\alpha)\cup \{*\}$ and say that $\alpha$ is a descendant of $\beta$, and $\beta$ is an ancestor of $\alpha$. We will sometimes denote the set of leaves $\Natural^r$ of $\A$ by $\LL(\A)$. For any $\alpha, \beta\in \A$, let
\begin{equation}
\alpha\wedge\beta 
:=
 |p(\alpha) \cap p(\beta)  |
\label{wedge}
\end{equation}
be the number of common vertices in the paths from the root to the vertices $\alpha$ and  $\beta$. In other words, $\alpha \wedge \beta$ is the distance of the lowest common ancestor of $\alpha$ and $\beta$ from the root. 

Let us now consider $r+1$ disjoint intervals $I_0,\ldots,I_r\subseteq[0,q^*]$. We allow all endpoint conventions, so each interval may have any of the forms
\begin{equation}
[a,b],\qquad [a,b),\qquad (a,b],\qquad (a,b),\qquad \{a\}.
\label{Ips}
\end{equation}
The intervals are ordered: if $p<p'$, then $x<y$ for every $x\in I_p$ and $y\in I_{p'}$. We assume that
\begin{equation}
{\rm supp}(\zeta)\subseteq \bigcup_{0\leq p\leq r} I_p \ \mbox{ and }\ \zeta(I_p) >0 \ \mbox{ for all }\ 0\leq p\leq r.
\label{intervals}
\end{equation}
Endpoints may be atoms of $\zeta$, and a singleton interval may consist of any atom, whether or not it is isolated. The intervals discretize the overlap values; when the overlap takes only finitely many values, we can discretize exactly by taking these values as singleton intervals.

Put $q_p=\inf I_p$. These infima are nondecreasing but need not be strictly increasing: for example, adjacent intervals $\{b\}$ and $(b,c)$ have the same infimum. To specify the cluster levels without ambiguity, define the bin index
\[
d_I(x)=p\qquad\mbox{for }x\in I_p.
\]
This function is nondecreasing on the union of the intervals, and we extend it to a nondecreasing $\{0,\ldots,r\}$-valued function on $[0,q^*]$. The extension across gaps does not affect the sampled overlaps, which belong to ${\rm supp}(\zeta)$ almost surely. For $0\leq p\leq r$, define the depth-$p$ cluster relation by
\begin{equation}
\sigma\sim_p^I\sigma'
\quad\Longleftrightarrow\quad
d_I(\sigma\cdot\sigma')\geq p.
\label{levelclusters}
\end{equation}
Since $d_I$ is nondecreasing, the discretized overlaps are ultrametric as well. Also, $d_I(q^*)=r$, so these are nested equivalence relations, with the depth-$0$ cluster equal to the entire support of $G$. For adjacent bins $\{b\}$ and $(b,c)$, the corresponding relations are $\sigma\cdot\sigma'\geq b$ and $\sigma\cdot\sigma'>b$, respectively. This notation includes both strict and non-strict overlap thresholds.

Later on we will need the sequence
\begin{equation}
0= \zeta_{-1} <\zeta_0 <\ldots < \zeta_{r-1} <\zeta_r = 1
\label{zetas}
\end{equation}
such that $\zeta_p - \zeta_{p-1} = \zeta(I_p)$ for $0\leq p\leq r.$ In particular,
\[
\zeta\{x:d_I(x)\geq p\}=1-\zeta_{p-1}.
\]
Let us now enumerate the clusters defined by (\ref{levelclusters}) according to Gibbs' weights as follows. Let $H_{*}$ be the entire support of $G$ so that $V_* = G(H_*) =1$. Next, the support is split into depth-$1$ clusters $(H_n)_{n\geq 1}$, which are then enumerated in the non-increasing order of their weights $V_n = G(H_n)$, 
\begin{equation}
V_1 \geq V_2 \geq \ldots \geq V_n \geq \ldots \geq 0.
\label{purelabelsfirst}
\end{equation}
We then continue recursively over $p\leq r-1$ and enumerate the depth-$(p+1)$ subclusters $(H_{\alpha n})_{n\geq 1}$ of a cluster $H_\alpha$ for $\alpha\in \Natural^p$ in the non-increasing order of their weights $V_{\alpha n} = G(H_{\alpha n})$, 
\begin{equation}
V_{\alpha 1} \geq V_{\alpha 2} \geq \ldots \geq V_{\alpha n} \geq \ldots \geq 0.
\label{purelabels}
\end{equation}
It is a well-known fact that each cluster $H_\alpha$ is split into infinitely many subclusters $(H_{\alpha n})_{n\geq 1}$ and their weights are all different and not equal to zero -- this is another consequence of the Ghirlanda-Guerra identities. Therefore, all the inequalities in (\ref{purelabelsfirst}) and (\ref{purelabels}) are strict. More specifically, it is well known that the Ghirlanda-Guerra identities imply that the cluster weights 
\begin{equation}
V = (V_\alpha)_{\alpha\in \A}
\label{Vs}
\end{equation}
are generated by the Ruelle probability cascades (RPC) \cite{Ruelle}, with parameters (\ref{zetas}). Indeed, the discretized overlaps $d_I(R_{\ell,\ell'})$ satisfy the Ghirlanda-Guerra identities, since $d_I$ is measurable, and their law depends only on the bin masses. This will be reviewed in Section \ref{Sec4label} (see also Chapter 2 in \cite{SKmodel}). We will call the depth-$r$ clusters $H_\alpha$ indexed by the leaves $\alpha\in \LL(\A) = \Natural^r$ the \emph{pure states}. Of course, if $\zeta(\{q^*\})>0$ then one can take $I_r = \{q^*\}$ in (\ref{Ips}) to ensure that the pure states are exactly the atoms of $G$. (For a way to construct pure states for the non-asymptotic Gibbs measure $G_N$ in (\ref{GibbsN}), see \cite{Tal-New}.)

Notice that the diameter of a pure state $H_\alpha$ for $\alpha\in\Natural^r$ can be bounded in $L_2$ by 
$$
{\rm diam}(H_\alpha) \leq \sqrt{2(q^*-q_r)},
$$
and when $q_r$ is close to $q^*$, these clusters are small and can be well approximated by one point, for example, the $G$-barycenter of the cluster. We can take these barycenters as an approximate definition of pure states but, in order not to lose any information, we will encode a pure state by an infinite sample as follows. First of all, notice that sampling from $G$ can now be done in two steps:
\begin{enumerate}
\item Choose $\alpha\in \LL(\A) = \Natural^r$ according to the weights $(V_\alpha)_{\alpha\in\Natural^r}$.
\item Sample from the pure state $H_\alpha$ according to the conditional distribution 
\begin{equation}
G_\alpha(\ \cdot\ ) = \frac{G( \ \cdot\ \cap H_\alpha)}{G(H_\alpha)}.
\end{equation}
\end{enumerate}
For each $\alpha\in \LL(\A) = \Natural^r$, let us consider an i.i.d. sample $(\sigma^{\alpha \ell})_{\ell \geq 1}$ with the distribution $G_\alpha$ and let  these samples be independent over such $\alpha$. As in (\ref{SpinsEll}), let us consider the sequence of spins 
\begin{equation}
S(\sigma^{\alpha \ell}) = \Bigl( 2 \Ind\Bigl(x_{i,\alpha \ell} \leq \frac{1+ \sigma^{\alpha \ell}(v_i) }{2}\Bigr) -1 \Bigr)_{i\geq 1}
\label{SpinsPure}
\end{equation}
generated by $\sigma^{\alpha \ell}$ and let
\begin{equation}
S_\alpha = (S(\sigma^{\alpha \ell}))_{\ell\geq 1}.
\label{Salpha}
\end{equation}
This array of spins completely encodes the pure state $H_\alpha$ for all practical purposes, if we remember that our main object of interest is the array of spins (\ref{sigma}) generated by the measure $G$.

To state our main result, it remains to recall the definition of hierarchical exchangeability introduced in \cite{AP}. Consider the following family of maps on the leaves $\Natural^r$ of the tree $\A$,
\begin{equation}
{\cal H} = \bigl\{ \pi: \Natural^r\to \Natural^r \,\bigr|\, \pi \mbox{ is a bijection}, \pi(\alpha)\wedge \pi(\beta) = \alpha\wedge\beta \mbox{ for all } \alpha,\beta\in \Natural^r \bigr\}.
\label{setH}
\end{equation}
As explained in \cite{AP}, the condition $\pi(\alpha)\wedge \pi(\beta) = \alpha\wedge\beta$ simply means that the genealogy on the tree is preserved after the permutation and such $\pi$ can be realized as a recursive rearrangement of children of each vertex starting from the root. We say that an array of random variables $(X_\alpha)_{\alpha\in \Natural^r}$ taking values in a standard Borel space is \emph{hierarchically exchangeable} if 
\begin{equation}
\bigl(X_{\pi(\alpha)} \bigr)_{\alpha\in \Natural^r}
\stackrel{d}{=}
\bigl(X_\alpha \bigr)_{\alpha\in \Natural^r}
\label{HexchDF}
\end{equation}
for all $\pi\in {\cal H}$. Our main result will be the following structure theorem for the Gibbs measure $G$. 
\begin{theorem}\label{Th1} If (\ref{GG}) holds then the array (\ref{Salpha}) of spins $(S_\alpha)_{\alpha\in \Natural^r}$ within pure states is hierarchically exchangeable and independent of the cluster weights $(V_\alpha)_{\alpha\in\A}$ in (\ref{Vs}).
\end{theorem}
If we write $S_\alpha = (S_{\alpha,i})_{i\geq 1}$, by making the dependence on the spin index $i$ in (\ref{SpinsPure}) explicit, then it is obvious that the distribution of the array $(S_{\alpha,i})$ is also invariant under the permutation of spins, 
\begin{equation}
\bigl(S_{\pi(\alpha), \rho(i)} \bigr)_{\alpha\in \Natural^r, i\in \Natural}
\stackrel{d}{=}
\bigl(S_{\alpha,i} \bigr)_{\alpha\in \Natural^r, i\in\Natural}
\label{HexchAH}
\end{equation}
for all $\pi\in {\cal H}$ and all bijections $\rho:\Natural \to \Natural$. The Aldous-Hoover representation was generalized to such hierarchically exchangeable arrays in \cite{AP} and, in particular, Theorem 2 in \cite{AP} implies the following.
\begin{corollary}
If (\ref{GG}) holds then the array $(S_{\alpha,i})_{\alpha\in \Natural^r, i\in \Natural}$ can be generated in distribution as
\begin{equation}
S_{\alpha,i} = f\bigl( \omega_*, (\omega_{\beta})_{\beta\in p(\alpha)}, \omega_*^i, (\omega_{\beta}^i)_{\beta\in p(\alpha)} \bigr),
\label{sigmaAH}
\end{equation}
where $f: [0,1]^{2(r+1)} \to\{-1,+1\}^\Natural$ is a measurable function and $\omega_\alpha,\omega_\alpha^i$ for $\alpha\in\A$ and $i\in \Natural$ are i.i.d. random variables with the uniform distribution on $[0,1]$. 
\end{corollary}
Note a slight difference in notation here and in \cite{AP} -- in this paper we chose not to include the root $*$ in the path (\ref{pathtoleaf}) while in \cite{AP} it was included. This is why we write $\omega_*$ and $\omega_*^i$ in (\ref{sigmaAH}) separately. Let us now explain the connection of the representation (\ref{sigmaAH}) to the M\'ezard-Parisi ansatz and what seems to be the main obstacle left. First of all, if we denote the barycenter of the pure state $H_\alpha$ by 
\begin{equation}
\bs^\alpha = \int_{H_\alpha} \sigma\, dG_\alpha(\sigma)
\end{equation}
then, by the strong law of large numbers, (\ref{SpinsPure}) implies that 
\begin{equation}
m^{\alpha} =(m_i^\alpha)_{i\geq 1}: = \bigl(\bs^\alpha(v_i)\bigr)_{i\geq 1} = \lim_{n\to\infty} \frac{1}{n} \sum_{\ell=1}^n S(\sigma^{\alpha \ell})
\label{magnet}
\end{equation}
almost surely. In the case when the pure state consists of one point $\bs^\alpha$ (for example, we mentioned above that if $\zeta(\{q^*\})>0$ and we choose $I_r = \{q^*\}$ then all pure states will be points) the vector $m^\alpha$ is called the \emph{magnetization} inside the pure state $\alpha$, otherwise, we can view it as an approximate notion of magnetization. The representation (\ref{sigmaAH}) and (\ref{magnet}) imply that
\begin{equation}
m^\alpha_i = m\bigl( \omega_*, (\omega_{\beta})_{\beta\in p(\alpha)}, \omega_*^i, (\omega_{\beta}^i)_{\beta\in p(\alpha)} \bigr)
\label{MAH}
\end{equation}
for some measurable function $m: [0,1]^{2(r+1)} \to [-1,1]$. 

What the M\'ezard-Parisi ansatz predicts is that, when $r$ is getting large and all the intervals $I_p$ in (\ref{Ips}) are getting small (which means that the $r$-step RSB scheme gives a good approximation of the overlap distribution), the magnetizations inside the pure states can be generated approximately (in the sense of distribution) by 
\begin{equation}
m^\alpha_i = m\bigl(\omega_*^i, (\omega_{\beta}^i)_{\beta\in p(\alpha)} \bigr)
\label{MAH2}
\end{equation}
for some measurable function $m: [0,1]^{r+1} \to [-1,1]$. As we already mentioned above, only the so called $1$-RSB case corresponding to $r=1$ was described in detail in \cite{Mezard} (see Section V there), but the general case is just a natural extension. The function $m$ is the order parameter of the M\'ezard-Parisi ansatz in the sense that one can express the free energy by some variational formula in terms of $m$. Obviously, (\ref{MAH2}) can hold only if the spin magnetizations are generated independently over the spin index $i\geq 1$ within pure states (which was, in fact, an assumption in \cite{Mezard}), but this assumption can be relaxed and the M\'ezard-Parisi formula for the free energy can be proved using the approach in Theorem $2$ in \cite{Pspins} under a slightly weaker hypothesis that the magnetizations inside the pure states are generated approximately by 
\begin{equation}
m^\alpha_i = m\bigl(\omega_*, \omega_*^i, (\omega_{\beta}^i)_{\beta\in p(\alpha)} \bigr)
\label{MAH3}
\end{equation}
for some measurable function $m: [0,1]^{r+2} \to [-1,1]$. The difference between (\ref{MAH}) and (\ref{MAH3}) can be informally expressed as follows. In (\ref{MAH3}), we have one (random) function $m(\omega_*, \ \cdot\ ,\ \cdot\ )$ that is used to generate spin magnetizations $m^\alpha_i$ in each pure state $\alpha$ using the randomness $\omega_*^i, (\omega_{\beta}^i)_{\beta\in p(\alpha)}$ along the path from the root to $\alpha$. In (\ref{MAH}), for each pure state $\alpha$ we first generate its own function $m(\omega_*, (\omega_{\beta})_{\beta\in p(\alpha)}, \ \cdot\ ,\ \cdot\ )$ in a hierarchically symmetric fashion and then use it to generate spin magnetizations inside that pure state. 

A possible way to go from (\ref{MAH}) to (\ref{MAH3}) is to show that \emph{multi-overlaps are functions of the overlaps}, which means the following. Let us consider $n$ pure state indices $\alpha_1,\ldots,\alpha_n \in \Natural^r$. If we compare the representations of $m^\alpha_i$ in terms of the barycenter $\bs^\alpha$ in (\ref{magnet}) and in terms of the function $m$ in (\ref{MAH}) then the so called \emph{multi-overlap} between these $n$ barycenters can be written as
$$
R_{\alpha_1,\ldots,\alpha_n} :=
\int \prod_{\ell\leq n}\bs^{\alpha_\ell}(v) \,dv
=
\e_i \prod_{\ell\leq n} m\bigl( \omega_*, (\omega_{\beta})_{\beta\in p(\alpha_\ell)}, \omega_*^i, (\omega_{\beta}^i)_{\beta\in p(\alpha_\ell)} \bigr),
$$
where $\e_i$ denotes the average in the random variables that depend on the spin index $i$. If (\ref{MAH3}) holds then, similarly,
$$
R_{\alpha_1,\ldots,\alpha_n} :=
\int \prod_{\ell\leq n}\bs^{\alpha_\ell}(v) \,dv
=
\e_i \prod_{\ell\leq n} m\bigl( \omega_*, \omega_*^i, (\omega_{\beta}^i)_{\beta\in p(\alpha_\ell)} \bigr),
$$
which clearly depends only on $(\alpha_\ell\wedge\alpha_{\ell'})_{1\leq \ell,\ell'\leq n}$. In the opposite direction, it is also not difficult to show that if $R_{\alpha_1,\ldots,\alpha_n}$ depends only on $(\alpha_\ell\wedge\alpha_{\ell'})_{1\leq \ell,\ell'\leq n}$ for all $n\geq 2$ then (\ref{MAH}) can be replaced by (\ref{MAH3}). Of course, in the $r$-step RSB approximation, $\alpha_\ell\wedge\alpha_{\ell'}$ describes the overlap $\bs^{\alpha_\ell}\cdot\bs^{\alpha_{\ell'}}$ only approximately, so the statement ``multi-overlaps are functions of overlaps" should be understood in an approximate sense for a finite $r$-step RSB approximation and should only become exact as $r$ goes to infinity, or if the distribution of the overlap is indeed concentrated on $r+1$ points. 

\smallskip
In the next section, we will begin with a review of the Ghirlanda-Guerra identities. In Section \ref{Sec3label}, we will prove some analogue of Theorem \ref{Th1} at the level of the sample from the Gibbs measure rather than working with the pure states directly. In Section \ref{Sec4label}, we will prove a technical result about the weights in the Ruelle probability cascades and, in Section \ref{Sec5label}, we will deduce Theorem \ref{Th1} from the main result in Section \ref{Sec3label} by sending the sample size to infinity.

\smallskip
\noindent
\textbf{Acknowledgements.} The author would like to thank the referees for very thorough reviews and many suggestions to improve the paper.

\section{The Ghirlanda-Guerra identities}

In this section, we will explain in what sense the Ghirlanda-Guerra identities are valid in diluted models, and we will use the example of the random $K$-sat model (\ref{Ham2}) for this purpose. For each $p\geq 1$, let us consider the process $g_p(\sigma)$ on $\Sigma_N = \{-1,+1\}^N$ given by
\begin{equation}
g_{p}(\sigma)
=
\frac{1}{N^{p/2}}
\sum_{i_1,\ldots,i_p = 1}^N g_{i_1,\ldots, i_p} \sigma_{i_1}\ldots\sigma_{i_p},
\label{mixedppert}
\end{equation}
where $(g_{i_1,\ldots, i_p})$ are i.i.d. standard Gaussian random variables, and define
\begin{equation}
g(\sigma) = \sum_{p\geq 1} 2^{-p} x_p\smsp g_{p}(\sigma)
\label{mixedHpert}
\end{equation}
for parameters $(x_p)_{p\geq 1}$ that take values in the interval $x_p\in[0,3]$ for all $p\geq 1$. It is easy to check that the variance of this Gaussian process satisfies $\e g(\sigma)^2 \leq 3.$ Given the Hamiltonian $H_N(\sigma)$ in (\ref{Ham2}), let us consider the perturbed Hamiltonian 
\begin{equation}
H_N^{\mathrm{pert}}(\sigma) = H_N(\sigma) - \frac{s}{\beta} g(\sigma)
\label{Hpert}
\end{equation} 
for some parameter $s\geq 0.$ It is easy to see, using Jensen's inequality on each side, that
\begin{align*}
\frac{1}{N}\smsp \e\log \sum_{\sigma\in\varSigma_N} \exp\bigl(-\beta  H_N(\sigma)\bigr)
& \leq \
\frac{1}{N}\smsp \e\log \sum_{\sigma\in\varSigma_N} \exp\bigl(-\beta H_N^{\mathrm{pert}}(\sigma) \bigr)
\\
& \leq \
\frac{1}{N}\smsp \e\log \sum_{\sigma\in\varSigma_N} \exp\bigl(-\beta  H_N(\sigma) \bigr)
+ \frac{3s^2}{2N}.
\end{align*}
Therefore, if we let $s$ in (\ref{Hpert}) depend on $N$, $s=s_N$, in such a way that
\begin{equation}
\lim_{N\to\infty} N^{-1} s_N^2 = 0, 
\label{tN}
\end{equation} 
then the limit of the free energy is not affected by the perturbation term $(s/\beta)g(\sigma).$ Since our ultimate goal is to find the formula for the free energy in the limit $N\to\infty$, adding a perturbation term is allowed if it helps us in some other way. Of course, the real purpose of adding the perturbation term is to obtain the Ghirlanda-Guerra identities for the Gibbs measure
\begin{equation}
G_N(\sigma) = \frac{\exp (-\beta  H_N^{\mathrm{pert}}(\sigma) )}{Z_N}
\,\mbox{ where }\,
Z_N = \sum_{\sigma\in\varSigma_{N}} \exp \bigl(-\beta   H_N^{\mathrm{pert}}(\sigma) \bigr),
\label{GNpert}
\end{equation}
which now corresponds to the perturbed Hamiltonian (\ref{Hpert}). In other words, even though the perturbation term does not affect the free energy, it will affect the Gibbs measure `in a good way' by regularizing it and forcing it to satisfy the Ghirlanda-Guerra identities. From now on, everything is defined with respect to this Gibbs measure corresponding to the perturbed Hamiltonian, such as any limit $\mu\in\M$ and the corresponding asymptotic Gibbs measure $G$.

Since we will soon pass to the limit $N\to\infty$, it should not cause any confusion if we temporarily denote by $\la\cdot\ra$ the average with respect to $G_N^{\otimes\infty}$ in (\ref{GNpert}), let $(\sigma^\ell)_{\ell\geq 1}$ be a sequence of replicas from $G_N$ and denote by 
\begin{equation}
R_{\ell,\ell'} = \frac{1}{N} \sum_{i=1}^N \sigma_i^\ell \sigma_i^{\ell'}
\label{overlapN}
\end{equation}
the overlap between replicas $\sigma^\ell$ and $\sigma^{\ell'}$. Let us consider the function
\begin{equation}
\varphi = \log Z_N =
\log \sum_{\sigma\in\varSigma_N} \exp \bigl( -\beta H_N(\sigma) + s g(\sigma)\bigr)
= \log \sum_{\sigma\in\varSigma_N} \exp \bigl( -\beta H_N^{\mathrm{pert}}(\sigma)\bigr),
\label{theta}
\end{equation}
viewed as a random function $\varphi = \varphi\bigl((x_p)\bigr)$ of the parameters $(x_p)$ in (\ref{mixedHpert}), and suppose that
\begin{equation}
\sup\Bigl\{ \e |\varphi - \e \varphi | \ \bigr| \ 0\leq x_p\leq 3, p\geq 1\Bigr\}\leq v_N(s)
\label{vt}
\end{equation}
for some function $v_N(s)$ that describes how well $\varphi((x_p))$ is concentrated around its expected value uniformly over all possible choices of the parameters $(x_p)$ from the interval $[0,3].$  Now, for any $n\geq 2, p\geq 1$ and any function $f=f(\sigma^1,\ldots,\sigma^n)$ on $\Sigma_N^n$ uniformly bounded by $1$, let us define
\begin{equation}
\varDelta(f,n,p) = 
\Bigl|
\e  \bigl\la f R_{1,n+1}^p \bigr\ra -  \frac{1}{n}\e \bigl\la f \bigr\ra \smsp \e\bigl\la R_{1,2}^p\bigr\ra - \frac{1}{n}\sum_{\ell=2}^{n}\e \bigl\la f R_{1,\ell}^p\bigr\ra
\Bigr|.
\label{GGfinite}
\end{equation}
Let us now think of $(x_p)_{p\geq 1}$ as a sequence of i.i.d. random variables with the uniform distribution on $[1,2]$ and denote by $\e_x$ the expectation with respect to such sequence. Here is one common  formulation of the Ghirlanda-Guerra identities from Theorem 3.2 in \cite{SKmodel}.
\begin{theorem}\label{ThGG} 
Suppose that the parameter $s$ in (\ref{Hpert}) depends on $N$,  $s=s_N$, and the sequence $(s_N)$ satisfies $\lim_{N\to\infty} s_N=\infty$ and $\lim_{N\to\infty} s_N^{-2} v_N(s_N) = 0$. Then
\begin{equation}
\lim_{N\to\infty} \e_x \smsp \varDelta(f,n,p) = 0
\label{GGxlim}
\end{equation} 
for any $p\geq 1, n\geq 2$ and any measurable function $f$ such that $\|f\|_\infty\leq 1$.
\end{theorem}
Of course, since the space $\Sigma_N$ changes with $N$, the function $f$ here is really a sequence $f=f_N$ such that $\|f_N\|_\infty\leq 1$ for all $N\geq 1$. 

We will show below that, in the setting of the $K$-sat model, one can find a sequence $(s_N)$ that satisfies (\ref{tN}) and the conditions in Theorem \ref{ThGG}. However, first let us recall how one can go from (\ref{GGxlim}) to (\ref{GG}) for any asymptotic Gibbs measure $G$ that arises in the limit (as explained in the introduction) from a sequence of measures $G_N$ that satisfy (\ref{GGxlim}). Simply, we consider the collection $\cal F$ of all triples $(f,n,p)$ such that $p\geq 1, n\geq 2$ and $f = \prod_{(i,\ell)\in F}\sigma_i^\ell$ for a finite subset $F\subseteq \Natural \times \{1,\ldots, n\}.$ This is a countable collection, so we can enumerate it, ${\cal F} = \{(f_j,n_j,p_j) \ |\ j\geq 1\}$, and consider
$$
\Delta_N(x) = \sum_{j\geq 1} 2^{-j} \Delta(f_j,n_j,p_j).
$$
Then (\ref{GGxlim}) implies that $\lim_{N\to\infty} \e_x \Delta_N(x) = 0$ and, as a consequence, we can choose a sequence $x^N = (x_p^N)_{p\geq 1}$ changing with $N$ such that $\lim_{N\to\infty} \Delta_N(x^N) = 0.$ Therefore, if we now define the perturbation (\ref{mixedHpert}) and the Gibbs measure (\ref{GNpert}) with this choice of parameters $x^N$ that depend on $N$, we get 
$$
\lim_{N\to\infty}
\Bigl|
\e  \bigl\la f R_{1,n+1}^p \bigr\ra -  \frac{1}{n}\e \bigl\la f \bigr\ra \smsp \e\bigl\la R_{1,2}^p\bigr\ra - \frac{1}{n}\sum_{\ell=2}^{n}\e \bigl\la f R_{1,\ell}^p\bigr\ra
\Bigr|
=0
$$
for any $(f,n,p)\in {\cal F}$. It should be obvious that this implies (\ref{GG}) for any asymptotic Gibbs measure $G$ corresponding to a limit $\mu\in\cal M$ of $(\mu_N)$ in (\ref{muN}) over any subsequence. The fact that the overlaps in (\ref{overlapN}) converge in distribution to the overlap in (\ref{overlapinfty}) over the same subsequence can be easily seen by computing their joint moments using the symmetry between sites (see the introduction in \cite{Pspins} for details). Moreover, the identities (\ref{GG}) for $\psi(x) = x^p$ and $f$ given by a product of finitely many spins, clearly, imply (\ref{GG}) for any $f$ and $\psi$. (Finally, let us point out that, even though the Ghirlanda-Guerra identities are typically proved via the above perturbation, in the mixed $p$-spin models they can be proved without any perturbation, see \cite{PGGmixed} or Section 3.7 in \cite{SKmodel}.)

Let us check the conditions of Theorem \ref{ThGG} in the random $K$-sat model.
\begin{lemma}
For the $K$-sat Hamiltonian (\ref{Ham2}), both (\ref{tN}) and the conditions in Theorem \ref{ThGG} are satisfied with $s_N = N^{\gamma}$ for any $\gamma\in (1/4, 1/2).$
\end{lemma}
\textbf{Proof.} We need to estimate the left hand side of (\ref{vt}) with $H_N(\sigma)$ given by (\ref{Ham2}). We will separate various sources of randomness as follows. For a function $\varphi =\varphi(X,Y)$ of two independent random variables $X$ and $Y$, by the triangle inequality and Jensen's inequality,
$$
\e |\varphi - \e \varphi| \leq \e |\varphi - \e_X \varphi| + \e | \e_X \varphi - \e \varphi| 
\leq \e |\varphi - \e_X \varphi| + \e |  \varphi - \e_Y \varphi|, 
$$ 
where $\e_X$ and $\e_Y$ denote the expectation in $X$ and $Y$ only. Similarly, for a function $\varphi =\varphi(X,Y,Z)$ of three independent random variables,
$$
\e |\varphi - \e \varphi| \leq  \e |  \varphi - \e_X \varphi| + \e |\varphi - \e_Y \varphi| + \e |  \varphi - \e_Z \varphi|. 
$$ 
In the case of the function (\ref{theta}), these three sources of randomness will come from the perturbation term $g(\sigma)$, the Poisson random variable $\pi(\alpha N)$, and the sequence of Rademacher random variables $(\eps_{j,k})$ and random indices $(i_{j,k})$. We will write the corresponding expectations by $\e_g$, $\e_\pi$ and $\e_\theta$ correspondingly, so that 
$$
\e |\varphi - \e \varphi| \leq  \e |  \varphi - \e_g \varphi| + \e |\varphi - \e_\pi \varphi| + \e |  \varphi - \e_\theta \varphi|. 
$$ 
In each term, we will first fix all other randomness and estimate $\e_g |  \varphi - \e_g \varphi|$, $\e_\pi |\varphi - \e_\pi \varphi|$  and  $\e_\theta |  \varphi - \e_\theta \varphi|$. The first one can be estimated using the standard Gaussian concentration (see e.g. Theorem 1.2 in \cite{SKmodel}). Since the variance of $sg(\sigma)$ is bounded by $3s^2$, we get $\e_g |  \varphi - \e_g \varphi| \leq L s$ for some absolute constant $L$. This gives $\e |  \varphi - \e_g \varphi| \leq Ls$. To estimate the last two terms, we will use the fact that each term in (\ref{Ham2}) for a fixed $k$,
\begin{equation}
\theta_k(\sigma) = \prod_{1\leq j\leq K} \frac{1+\eps_{j,k} \sigma_{i_{j,k}}}{2},
\label{thetak}
\end{equation}
is bounded uniformly by $1$. First of all, if $\pi_1$ and $\pi_2$ are two independent copies of $\pi(\alpha N)$, and we think of $\varphi$ for a moment as a function $\varphi(\pi(\alpha N))$ of $\pi(\alpha N)$ only, then
$$
\e_\pi |\varphi - \e_\pi \varphi| \leq \e_\pi |\varphi(\pi_1) - \varphi(\pi_2)| 
\leq \beta \e |\pi_1 - \pi_2| \leq 2\beta\sqrt{\alpha N}.
$$
This gives $\e |\varphi - \e_\pi \varphi| \leq 2\beta\sqrt{\alpha N}$. Finally, to estimate $\e_\theta |  \varphi - \e_\theta \varphi|$, we can use the standard martingale difference representation for $\varphi - \e_\theta \varphi = \sum_{k\leq \pi(\alpha N)} d_k$ by adding the randomness of one term (\ref{thetak}) at a time to obtain
$$
\e_\theta (\varphi - \e_\theta \varphi)^2 = \sum_{k\leq \pi(\alpha N)} \e_\theta d_k^2 \leq 4\beta^2 \pi(\alpha N).
$$
Therefore, $\e (\varphi - \e_\theta \varphi)^2 \leq 4\beta^2 \alpha N$ and $\e |  \varphi - \e_\theta \varphi| \leq 2\beta\sqrt{\alpha N}$. Combining all three estimates, we proved that $\e |\varphi - \e \varphi| \leq L s + 4\beta \sqrt{\alpha N}.$ Now it is easy to see that we can take $s_N = N^{\gamma}$ for any $\gamma\in (1/4, 1/2)$ to satisfy (\ref{tN}) and the conditions in Theorem \ref{ThGG}. 
\qed

\medskip
\medskip
\noindent
We now go back to the notations in the setting of asymptotic Gibbs measures in the introduction, and will end this section with the invariance property that will be the main tool in the proof of Theorem \ref{Th1}. Given $n\geq 1$, consider $n$ bounded measurable functions $f_1,\ldots, f_n: \Reals\to\Reals$ and define
 \begin{equation}
F(\sigma,\sigma^1,\ldots,\sigma^n) = f_1(\sigma\cdot\sigma^1)+\ldots+f_n(\sigma\cdot\sigma^n).
\label{F1}
\end{equation}
For $1\leq \ell\leq n$ we define
\begin{equation}
F_\ell(\sigma,\sigma^1,\ldots,\sigma^n) = F(\sigma,\sigma^1,\ldots,\sigma^n)
- f_\ell(\sigma\cdot\sigma^\ell)+ \e\bigl\la  f_\ell(R_{1,2}) \bigr\ra
\label{F2}
\end{equation}
Consider a finite index set $\T.$ Given a realization of the random measure $G$ and a sample $\sigma^1,\ldots,\sigma^n$ from $G$ let $(B_t)_{t\in\T}$ be a partition of the support of $G$ such that, for each $t\in\T$, the indicator $\Ind(\sigma\in B_t)$ is a measurable function of $(\sigma^\ell\cdot \sigma^{\ell'})_{\ell,\ell'\leq n}$ and $(\sigma\cdot\sigma^\ell)_{\ell\leq n}$. Let
\begin{equation}
\delta_t=\delta_t(\sigma^1,\ldots,\sigma^n)=G(B_t).
\label{WA}
\end{equation}
Let us define the map $T$ by
\begin{equation}
\delta=(\delta_t)_{t\in\T}\to T(\delta) = 
\Bigl(\frac{\la \Ind(\sigma\in B_t) \exp F(\sigma,\sigma^1,\ldots,\sigma^n )\ra_{\mathunderscore}}
{\la \exp F(\sigma,\sigma^1,\ldots,\sigma^n )\ra_{\mathunderscore}} \Bigr)_{t\in\T},
\label{TA}
\end{equation}
where $\la\cdot\ra_{\mathunderscore}$ denotes the average  with respect to the measure $G$ in $\sigma$ only for fixed  $\sigma^1,\ldots, \sigma^n$. The following result was proved in \cite{PUltra} (see also Theorem 2.19 in \cite{SKmodel}) as a consequence of the Ghirlanda-Guerra identities (\ref{GG}). Recall the definition of $S^n$ in (\ref{Sn}).
\begin{theorem}\label{Th2}
If (\ref{GG}) holds then, for any bounded measurable function $\Phi =\Phi(S^n, \delta)$,
\begin{equation}
\e\bigl \la  \Phi(S^n, \delta) \bigr\ra
=
\e\Bigl\la
\frac{ \Phi(S^n,T(\delta)) \exp \sum_{\ell=1}^{n} F_\ell(\sigma^\ell,\sigma^1,\ldots,\sigma^n)}
{\la\exp F(\sigma,\sigma^1,\ldots,\sigma^n)\ra_{\mathunderscore}^n}
\Bigr\ra.
\label{nA}
\end{equation}
\end{theorem}
This theorem was proved in \cite{PUltra} for the function $\Phi$ of the overlaps $(R_{\ell,\ell'})_{\ell,\ell'\leq n}$ instead of all spins $S^n$. This is because the Ghirlanda-Guerra identities in \cite{PUltra} were stated only for the function of the overlaps, while here we wrote them in (\ref{GG}) for a function of all spins. Otherwise, the proof of Theorem \ref{Th2} from (\ref{GG}) is identical to the one in \cite{PUltra}.

\section{At the level of replicas} \label{Sec3label}

The main work will be to prove some analogue of Theorem \ref{Th1} at the level of the replicas $\sigma^1,\ldots, \sigma^n$ sampled from the Gibbs measure $G$ described at the end of previous section, which will then imply Theorem \ref{Th1} by passing $n$ to infinity. Until further notice, however, $n$ will be fixed.

Let $\T$ be a finite rooted labelled tree of depth $r$. We will label the vertices of $\T$ by a finite subset of $\A$ in (\ref{Atree}) as follows. The root will again be labelled by $*$. Then, recursively for $p\leq r-1$, if a vertex at the distance $p$ from the root labelled by $t\in \Natural^p$ has $k_t$ children then we label them by $t 1,\ldots, t k_t \in \Natural^{p+1}$ (recall that for simplicity we write $tk$ for $(t,k)$). We identify the tree $\T$ with the set of vertex labels and use the same notation, $|t|, t\wedge s$, $t\succ s$ for $t,s\in \T$, as for the tree $\A$. We will denote by $\LL(\T)$ the set of leaves of $\T$ and consider a function
\begin{equation}
\PP: \{1,\ldots, n\} \to \LL(\T).
\end{equation}
We will call the pair $\CC=(\T,\PP)$ \emph{a configuration} if $\PP^{-1}(t) \not = \emptyset \,\mbox{ for all }\, t\in \LL(\T),$ i.e. at least one replica index is mapped into each leaf. Of course, this means that the cardinality $|\LL(\T)| \leq n$. The role of the function $\PP$ is to partition replica indices among the leaves of $\T$ and then use the tree structure to describe how replicas $\sigma^1,\ldots, \sigma^n$  cluster according to the overlap equivalence relations (\ref{levelclusters}) along the tree $\T$. More precisely, we will consider the event
\begin{equation}
\OC = \Bigl\{ (\sigma^1,\ldots,\sigma^n) \ \bigr|\ \sigma^\ell\cdot\sigma^{\ell'}\in I_{\PP(\ell)\wedge \PP(\ell')}  \mbox{ for all } 1\leq \ell,\ell' \leq n\Bigr\}.
\label{OCset}
\end{equation}
This event depends on the tree $\T$ via $\PP(\ell)\wedge \PP(\ell')$ and $I_{\PP(\ell)\wedge \PP(\ell')}$ is one of the intervals in (\ref{intervals}). In other words, on this event the overlap of replicas ``assigned by $\PP$" to the leaves $t,t'\in\LL(\T)$ is determined by the depth $t\wedge t'$ of their lowest common ancestor.

Let us assume from now on that the sample belongs to the event $\OC$. Then, we can use ultrametricity of the support of the measure $G$ to partition it in a natural way `along the tree $\T$' according to the overlaps with the replicas $\sigma^1,\ldots, \sigma^n$. For each $t\in \T$, let
\begin{equation}
\RR(t) = \bigl\{ 1\leq \ell \leq n \ |\  \PP(\ell) \succ t  \bigr\}
\label{Rt}
\end{equation}
be the set of replica indices assigned to the leaves which are descendants of $t$. Consider the sets
\begin{equation}
C_t = \bigl\{\sigma \ | \ d_I(\sigma\cdot \sigma^{\ell}) \geq |t|  \mbox{ for all } \ell\in \RR(t) \bigr\}.
\label{Ctall}
\end{equation}  
Since, obviously, $t'\wedge t'' \geq |t|$ for any $t',t''\succ t$, we have $d_I(\sigma^\ell\cdot\sigma^{\ell'}) \geq |t|$ for all $\ell,\ell' \in \RR(t)$ on the event $\OC$. By ultrametricity, this implies that we can also write the set (\ref{Ctall}) as
\begin{equation}
C_t = \bigl\{\sigma \ | \ d_I(\sigma\cdot \sigma^{\ell}) \geq |t|  \mbox{ for any } \ell\in \RR(t) \bigr\}.
\label{Ctany}
\end{equation}  
This makes it obvious that the sets $C_t$ are nested, $C_{t'} \subseteq C_t$ for $t'\succ t$. Another simple property is that the sets indexed by the children of $t$ are disjoint subsets of $C_t$,
\begin{equation}
C_{tk} \cap C_{tk'} = \emptyset \ \mbox{ for all }\ k\not = k'\leq k_t
\label{Cempty}
\end{equation}
(recall that $k_t$ is the number of children of $t\in\T$). To see this, if we take $\ell \in \RR(tk)$ and $\ell'\in \RR(tk')$ then $\sigma^\ell\cdot\sigma^{\ell'} \in I_{|t|}$, so $d_I(\sigma^\ell\cdot\sigma^{\ell'})=|t|$, by (\ref{OCset}). On the other hand,
$$
\mbox{
$d_I(\sigma\cdot \sigma^{\ell}) \geq |t|+1$ for $\sigma\in C_{tk}$ and $d_I(\sigma\cdot \sigma^{\ell'}) \geq |t|+1$ for $\sigma\in C_{tk'},$
}
$$
so (\ref{Cempty}) again follows by ultrametricity. Let us now consider the sets $B_t := C_t$ for $t\in \LL(\T)$ and
\begin{equation}
B_t : = C_t \setminus \cup_{k\leq k_t} C_{tk} = \bigl\{\sigma \ |\ \sigma\cdot\sigma^\ell \in I_{|t|} \mbox{ for all } \ell\in \RR(t)\bigr\}
\label{Bt}
\end{equation}
for $t\in \T\setminus \LL(\T).$ On the event $\OC$, the collection $(B_t)_{t\in\T}$ forms a random partition of the support of the Gibbs measures $G$ and, by definition, the indicator $\Ind(\sigma\in B_t)$ depends only on the overlaps $(\sigma\cdot\sigma^\ell)_{\ell\leq n}$. Below, this will allow us to apply Theorem \ref{Th2} to this partition with some specific choice of function $f_1,\ldots, f_n$ in (\ref{F1}). 

Let us denote the Gibbs weights of the above sets by 
\begin{equation}
W_t = G(C_t) \ \mbox{ and } \
\delta_t = G(B_t) = W_t - \sum_{k\leq k_t} W_{tk}. 
\label{WC}
\end{equation}
It is obvious that two different configurations $\CC=(\T,\PP)$ and $\CC'=(\T',\PP')$ can result in the same event, $\OC={\cal O}({\cal C}')$, if we simply reshuffle the labels of $\T$ in a hierarchical way and then redefine $\PP$ accordingly. Later on, we will need to fix a special configuration among these, and this will be done using the cluster weights $W_t$ around the sample points, as follows. Consider the event 
\begin{equation}
\WC = \Bigl\{ (\sigma^1,\ldots,\sigma^n) \ \bigr|\ 
W_{t1}> \ldots > W_{tk_t} \mbox{ for all } t\in\T\setminus \LL(\T)
\Bigr\}.
\label{WCset}
\end{equation}
It is obvious that such ordering of the weights makes the events $\OC\cap \WC$ disjoint for different configurations $\CC$, and each sample $(\sigma^1,\ldots,\sigma^n)$ belongs to one and only one of these events. We will denote the corresponding configuration by $\CC_n = (\T_n,\PP_n)$,
\begin{equation}
\CC_n = \CC \Longleftrightarrow (\sigma^1,\ldots,\sigma^n) \in \OC\cap \WC,
\label{sampleconfig}
\end{equation} 
and call $\CC_n = (\T_n,\PP_n)$ the \emph{sample configuration}. In other words, $\CC_n$ is a function of $\sigma^1,\ldots,\sigma^n$ such that the tree is defined according to the overlap structure $(\sigma^\ell\cdot\sigma^{\ell'})_{\ell,\ell'\leq n}$ with the vertices labelled according to the weights of the neighborhoods of these replicas. The event $\WC$ and the sample configuration $\CC_n$ will not be used in this section, but will play an important role in the last section where they will be utilized to partition an event into disjoint events indexed by configurations $\CC$. 

For the remainder of this section, we will fix a configuration $\CC$ once and for all and, for simplicity of notation, will omit the dependence of $\OC$ on $\CC$ and write $\OO$ instead. Let us denote $\p(\ \cdot \ ) = \e\la\Ind(\ \cdot\ )\ra$ and let  
\begin{equation}
\p_\OO(\ \cdot \ ) = \frac{\p(\ \cdot\ \cap \OO)}{\p(\OO)}
\label{PO}
\end{equation}
be the conditional distribution given the event $\OO$. Since $n$ is fixed in this section, we will write $S$ to denote $S^n$ in (\ref{Sn}). Let 
\begin{equation}
\T_* : = \T\setminus \{*\} \ \mbox{ and }\  W = (W_t)_{t\in \T_*}.
\label{WCT}
\end{equation}
We exclude the root, because $W_*=1$. Theorem \ref{Th1} will follow from the main result of this section. 
\begin{theorem}\label{Th3} For any measurable sets $A$ and $B$,
\begin{equation}
\p_\OO( S \in A, W\in B) = \p_\OO( S \in A)\hspace{0.2mm} \p_\OO(W\in B).
\label{SWind}
\end{equation}
\end{theorem}
Since the weights $(W_t)$ and $(\delta_t)$ in (\ref{WC}) are functions of each other, the independence of $S$ and $W$ in (\ref{SWind}) is equivalent to independence of $S$ and $\delta$,
\begin{equation}
\p_\OO( S \in A, \delta\in B) = \p_\OO( S \in A)\hspace{0.2mm} \p_\OO(\delta\in B),
\label{SDind}
\end{equation}
where $\delta = (\delta_t)_{t\in\T_*}$. Again, we can exclude the root, because $\delta_* = 1-\sum_{t\in\T_*}\delta_t.$ The vector $\delta$ takes values in the open subset
\begin{equation}
\DD = \Bigl\{ (x_t)_{t\in \T_*} \ \bigr|\  \sum_{t\in \T_*} x_t<1 \mbox{ and all } x_t>0 \Bigr\}
\label{DD}
\end{equation}
of $\Reals^{|\T_*|}.$ Given a vector $a=(a_t)_{t\in \T_*}\in \Reals^{|\T_*|}$, let us define the map $T_a: \DD\to \DD$ by 
\begin{equation}
T_a(x) = \Bigl(\frac{x_t e^{a_t}}{\Delta_a(x)}\Bigr)_{t\in \T_*}
\ \mbox{ where } \
\Delta_a(x) = \sum_{t\in \T_*} x_t e^{a_t} + 1- \sum_{t\in \T_*} x_t.
\label{Tax}
\end{equation}
One can easily check that for $a,b\in \Reals^{|\T_*|}$ we have $T_a\circ T_b = T_{a+b}$ and, therefore, $T_a^{-1} = T_{-a}$. It is also easy to check that 
\begin{equation}
\Delta_a(T_{-a}(x)) = \frac{1}{\Delta_{-a}(x)}.
\label{DeltaT}
\end{equation}
Let us denote by $B_\eps(x)$ the open ball of radius $\eps$ in $ \Reals^{|\T_*|}$ centered at $x.$ Then the following holds.
\begin{lemma}\label{ThBeps}
For any $a=(a_t)_{t\in \T_*}\in \Reals^{|\T_*|}$ and $x\in \DD$,
\begin{equation}
\lim_{\eps\downarrow 0}  
\frac{\p_\OO(S\in A, \delta\in B_\eps(x))}{ \p_\OO(\delta\in B_\eps(x))}
=
\lim_{\eps\downarrow 0}
 \frac{\p_\OO(S\in A, T_a(\delta) \in B_\eps(x))}{ \p_\OO(T_a(\delta) \in B_\eps(x))}
\label{Beps}
\end{equation}
whenever either of the limits exists.
\end{lemma}\noindent
\textbf{Proof.}  As we mentioned above, we will apply Theorem \ref{Th2} to the partition $(B_t)_{t\in\T}$ in (\ref{Bt}) with the following choice of function $f_1,\ldots, f_n$ in (\ref{F1}). Let us consider an arbitrary function 
\begin{equation}
\ell: \T \to \{1,\ldots, n\}
\label{pickleaf}
\end{equation}
such that $\ell(t) \in \RR(t)$ in (\ref{Rt}) for all $t\in \T$. In other words, we pick one replica index $\ell(t)$ assigned to one of the leaves that are descendants of $t$. Consider a vector $b=(b_t)_{t\in \T}\in \Reals^{|\T|}$. For each replica index $1\leq \ell\leq n$, let 
\begin{equation}
\T_\ell = \bigl\{t\in \T \ |\  \ell(t) = \ell \bigr\}
\ \mbox{ and }\
f_\ell(x) = \sum_{t\in \T_\ell} b_t \Ind(x\in I_{|t|}).
\end{equation}
Then the function $F$ in (\ref{F1}) can be written as
$$
F(\sigma,\sigma^1,\ldots,\sigma^n) 
= \sum_{\ell \leq n} \sum_{t\in \T_\ell} b_t \Ind(\sigma\cdot\sigma^\ell  \in I_{|t|})
= \sum_{t\in \T} b_t \Ind(\sigma\cdot\sigma^{\ell(t)}  \in I_{|t|}).
$$
Let us fix $u\in \T$ and compute $\la \Ind(\sigma\in B_u) \exp F(\sigma,\sigma^1,\ldots,\sigma^n )\ra_{\mathunderscore}$. We will now fix $\sigma\in B_u$ and consider several different cases when $t$ belongs to different subsets of the tree $\T$. 
\begin{enumerate}
\item
First of all, if $t=u$ then $\Ind(\sigma\cdot\sigma^{\ell(t)}  \in I_{|t|}) = 1$ by the definition of $B_u$ in (\ref{Bt}).

\item
If $t\succ u$, $t\not =u,$ then $\ell(t) \in \RR(u)$ and $\sigma\cdot\sigma^{\ell(t)} \in I_{|u|}$, which implies that $\Ind(\sigma\cdot\sigma^{\ell(t)}  \in I_{|t|}) = 0$. 

\item
If neither $t\succ u$ nor $u\succ t$ then (on the event $\OO$) $\sigma^{\ell(t)}\cdot \sigma^{\ell(u)} \in I_{t\wedge u}$ and $t\wedge u < \min(|t|,|u|)$. Since for $\sigma\in B_u$ we have $\sigma\cdot\sigma^{\ell(u)}  \in I_{|u|}$, by ultrametricity, $\Ind(\sigma\cdot\sigma^{\ell(t)}  \in I_{|t|}) = 0.$

\item 
If $u\succ t$, $t\not = u$, then, in general, the answer depends on the choice of the function (\ref{pickleaf}) or, more specifically, on whether $\PP(\ell(u))\wedge \PP(\ell(t)) = |t|$ or $\PP(\ell(u))\wedge \PP(\ell(t)) > |t|$. In the first case, on the event $\OO$, $\sigma^{\ell(t)}\cdot \sigma^{\ell(u)} \in I_{|t|}$. Since for $\sigma\in B_u$ we have $\sigma\cdot\sigma^{\ell(u)}  \in I_{|u|}$ and $I_{|u|}$ lies strictly above $I_{|t|}$, by ultrametricity, $\sigma \cdot \sigma^{\ell(t)} = \sigma^{\ell(t)}\cdot \sigma^{\ell(u)} \in I_{|t|}$ and $\Ind(\sigma\cdot\sigma^{\ell(t)}  \in I_{|t|}) = 1$. In the second case, on the event $\OO$, $\sigma^{\ell(t)}\cdot \sigma^{\ell(u)}$ also lies strictly above $I_{|t|}$ and, therefore, by ultrametricity, $\sigma\cdot\sigma^{\ell(t)}$ lies strictly above $I_{|t|}$ and $\Ind(\sigma\cdot\sigma^{\ell(t)}  \in I_{|t|}) = 0$.
\end{enumerate}
Therefore, if we consider the set
$$
\T(u) = \Bigl\{ t\in \T   \, \bigr| \ u\succ t, t\not =u \, \mbox{ and }  \PP(\ell(u))\wedge \PP(\ell(t)) = |t| \Bigr\}
$$
then, for $\sigma\in B_u$ we have $F(\sigma,\sigma^1,\ldots,\sigma^n) = b_u + \sum_{t\in \T(u)} b_t$ and
$$
\bigl\la \Ind(\sigma\in B_u) \exp F(\sigma,\sigma^1,\ldots,\sigma^n )\bigr\ra_{\mathunderscore}
=
G(B_u) \exp\Bigl(b_u+   \sum_{t\in \T(u)} b_t \Bigr).
$$
Let us now set $b_{*}=0$ and recursively set $b_u = a_u -  \sum_{t\in \T(u)} b_t$ for $u\in \T_{*}$. Then,
\begin{align*}
\bigl\la \Ind(\sigma\in B_*) \exp F(\sigma,\sigma^1,\ldots,\sigma^n )\bigr\ra_{\mathunderscore}
&=\
\delta_* = 1-\sum_{t\in \T_*} \delta_t,
\\
\bigl\la \Ind(\sigma\in B_u) \exp F(\sigma,\sigma^1,\ldots,\sigma^n )\bigr\ra_{\mathunderscore}
&=\
\delta_u e^{a_u}
\ \mbox{ for }\ u\in \T_*.
\end{align*} 
Adding them up, we get
$$
\bigl\la \exp F(\sigma,\sigma^1,\ldots,\sigma^n )\bigr\ra_{\mathunderscore}
= \sum_{t\in\T_*}\delta_u e^{a_u} + 1-\sum_{t\in \T_*} \delta_t = \Delta_a(\delta).
$$
We showed that, with this choice of functions $f_1,\ldots, f_n$, the map $T$ in (\ref{TA}) coincides with the map $T_a$ in (\ref{Tax}) on the coordinates indexed by $t\in\T_*.$ Also, it is clear that, on the event $\OO$, the sum
$$
\sum_{\ell=1}^n F_\ell(\sigma^\ell,\sigma^1,\ldots,\sigma^n )
$$
is a constant, which we will denote by $\gamma(a)$. If we denote $Z_a(\delta) = e^{\gamma(a)}/\Delta_a(\delta)^n$ then Theorem \ref{Th2} implies that
$$
\e\bigl \la  \Ind\bigl(S\in A, \delta\in B_\eps(x) \bigr) \Ind_{\OO} \bigr\ra
=
\e\bigl\la
\Ind\bigl(S\in A,T_a(\delta)\in B_\eps(x)\bigr) Z_a(\delta) \Ind_{\OO} 
\bigr\ra.
$$
The same equality, obviously, holds without the event $\{S\in A\}$, which proves that
\begin{equation}
\frac{\e \la  \Ind(S\in A, \delta\in B_\eps(x)) \Ind_{\OO} \ra}{\e \la  \Ind(\delta\in B_\eps(x)) \Ind_{\OO} \ra}
=
\frac{\e\la \Ind(S\in A,T_a(\delta)\in B_\eps(x)) Z_a(\delta) \Ind_{\OO} \ra}{\e\la \Ind(T_a(\delta)\in B_\eps(x)) Z_a(\delta) \Ind_{\OO} \ra},
\label{ratioZ}
\end{equation}
if the denominator is not zero. When $T_a(\delta)\in B_\eps(x)$, by (\ref{DeltaT}),
$$
\frac{1}{\Delta_a(\delta)} = \Delta_{-a}(T_a(\delta)) \in \Delta_{-a}(B_\eps(x))
$$
and, hence, $Z_a(\delta) \in e^{\gamma(a)} \Delta_{-a}(B_\eps(x))^n$. As a result, as $\eps\downarrow 0$, the factor $Z_a(\delta)$ converges uniformly to a constant $e^{\gamma(a)} \Delta_{-a}(x)^n$ that will cancel out on the right hand side of (\ref{ratioZ}), yielding (\ref{Beps}).
\qed

\medskip
\medskip
\noindent
We will need one more technical result that will be postponed until the next section.
\begin{lemma}\label{LemAC}
The distribution $\p_\OO(\delta \in \ \cdot \ )$ of weights $\delta = (\delta_t)_{t\in\T_*}$ is absolutely continuous with respect to the Lebesgue measure on $\Reals^{|\T_*|}$.
\end{lemma}
We are now ready to prove Theorem \ref{Th3}.

\medskip
\noindent
\textbf{Proof of Theorem \ref{Th3}.} 
Let $p(x)$ be the Lebesgue density of the distribution $\p_\OO(\delta \in \, \cdot \, )$ and let $p_A(\delta)$ be the conditional expectation of the indicator $\Ind(S\in A)$ given $\delta$ under the measure $\p_\OO.$ Then, 
\begin{equation}
\p_\OO(S\in A, \delta\in B) = \int_B p_A(x) p(x) \, dx
\ \mbox{ and } \
\p_\OO(\delta\in B) = \int_B p(x) \,dx.
\label{pAp}
\end{equation}
To prove (\ref{SDind}), it is enough to show that $p_A(x)$ is a constant a.e. on the set $\{x: p(x)>0\}$. By the Lebesgue differentiation theorem (see Corollary 1.6 in \cite{Stein}), for almost every $x'\in\Reals^{|\T_*|}$,
\begin{align}
&
\lim_{\eps\downarrow 0} \,\frac{1}{|B_\eps(x')|}\int_{B_\eps(x')} \bigl| p_A(x)p(x) - p_A(x')p(x')\bigr| \,dx = 0,
\label{Leb1}
\\
&
\lim_{\eps\downarrow 0}\, \frac{1}{|B_\eps(x')|}\int_{B_\eps(x')} \bigl|p(x) - p(x') \bigr| \,dx = 0.
\label{Leb2}
\end{align}
If $p_A(x)$ is not a constant a.e. on $\{p(x)>0\}$ then we can find two points $x',x''$ for which both (\ref{Leb1}) and (\ref{Leb2}) hold and such that $p(x'), p(x'')>0$ and $p_A(x')\not = p_A(x'').$ We can also assume that $x',x''\in \DD$ in (\ref{DD}) since $\p_\OO(\delta\not\in \DD) = 0.$ First of all, equations (\ref{pAp}) -- (\ref{Leb2}) imply that the left hand side of (\ref{Beps}) is equal to
\begin{equation}
\lim_{\eps\downarrow 0}  \,
\frac{\p_\OO(S\in A, \delta\in B_\eps(x'))}{ \p_\OO(\delta\in B_\eps(x'))}
= p_A(x').
\label{BepsLHS}
\end{equation}
It is easy to check that if we take
$$
a_t = \log \frac{x_t'}{x_t''} - \log \frac{1-\sum_{t\in \T_*} x_t'}{1-\sum_{t\in \T_*} x_t''}
$$
for $t\in \T_*$ then $T_a(x'') = x'$ for $T_a$ defined in (\ref{Tax}). Equations (\ref{Beps}) and (\ref{BepsLHS}) imply that
\begin{equation}
\lim_{\eps\downarrow 0} \,  \frac{\p_\OO(S\in A, \delta\in T_{-a} (B_\eps(x')))}{ \p_\OO(\delta\in T_{-a} (B_\eps(x')))} = p_A(x').
\label{contr1}
\end{equation}
To finish the proof, we will follow the argument of Corollary 1.7 in \cite{Stein} and use the fact that the sets $T_{-a} (B_\eps(x'))$ are of bounded eccentricity. Namely, since all partial derivatives of $T_a$ are uniformly bounded in a small neighborhood of $x''$ and all partial derivatives of $T_a^{-1}=T_{-a}$ are uniformly bounded in a small neighborhood of $x'$, there exist some constants $c, C>0$ such that $B_{c \eps}(x'') \subseteq T_{-a} (B_\eps(x')) \subseteq B_{C \eps}(x'')$ for small $\eps>0$. Therefore,
\begin{align*}
\frac{1}{| T_{-a} (B_\eps(x')) |}\int_{ T_{-a} (B_\eps(x'))} \bigl|p(x) - p(x'')\bigr| \,dx\
& \leq\
\frac{1}{ |B_{c\eps}(x'')|}\int_{B_{C\eps}(x'')} \bigl|p(x) - p(x'')\bigr|\, dx 
\\
& =\ \frac{(C/c)^{|\T_*|}}{|B_{C\eps}(x'')|}\int_{B_{C\eps}(x'')} \bigl|p(x) - p(x'')\bigr| \,dx,
\end{align*}
and, using that (\ref{Leb2}) holds with $x''$ instead of $x'$, we get
$$
\lim_{\eps\downarrow 0} \, \frac{1}{| T_{-a} (B_\eps(x')) |}\int_{ T_{-a} (B_\eps(x'))} \bigl|p(x) - p(x'')\bigr| \,dx = 0.
$$
Similarly, using (\ref{Leb1}) with $x''$ instead of $x'$, 
$$
\lim_{\eps\downarrow 0} \, \frac{1}{| T_{-a} (B_\eps(x')) |}\int_{ T_{-a} (B_\eps(x'))} \bigl|p_A(x) p(x) - p_A(x'') p(x'')\bigr|\,dx = 0.
$$
These equations together with (\ref{pAp}) for $B = T_{-a} (B_\eps(x'))$ imply that
$$
\lim_{\eps\downarrow 0} \,  \frac{\p_\OO(S\in A, \delta\in T_{-a} (B_\eps(x')))}{ \p_\OO(\delta\in T_{-a}( B_\eps(x')))} = p_A(x'').
$$
Recalling (\ref{contr1}), we arrive at contradiction, $p_A(x') = p_A(x'')$.
\qed

\section{Absolute continuity of cluster weight distribution}\label{Sec4label}

In this section, we will prove Lemma \ref{LemAC}. First, we will reduce the problem to proving absolute continuity for the distribution of finitely many cluster weights $V_\alpha$ in (\ref{Vs}). Then we will recall that these cluster weights are generated by the Ruelle probability cascades thanks to the Ghirlanda-Guerra identities, so the proof of absolute continuity will be based solely on the properties of the RPC.  

Let $\CC=(\T,\PP)$ be a fixed configuration as in the previous section. With probability one, the vector of weights $W=(W_t)_{t\in\T_*}$ defined in (\ref{WC}) belongs to the open subset
\begin{equation}
\WW = \Bigl\{ (y_t)_{t\in \T_*} \ \bigr|\  \sum_{k\leq k_t} y_{tk}<y_t \mbox{ for } t\in\T\setminus\LL(\T) \mbox{ and all } y_t>0 \Bigr\}
\label{WW}
\end{equation}
of $\Reals^{|\T_*|},$ where we set $y_* = 1.$ The map given by $x_t = y_t -\sum_{k\leq k_t} y_{tk}$ for $t\in \T_*$ is a linear bijection between $\WW$ and the set $\DD$ defined in (\ref{DD}). Recall that this is precisely the relationship between the weights $W=(W_t)_{t\in\T_*}$ and $\delta= (\delta_t)_{t\in\T_*}$ in (\ref{WC}). Therefore, in order to prove Lemma \ref{LemAC}, it is enough to prove that the distribution $\p_\OO(W\in\ \cdot\ )$ is absolutely continuous with respect to the Lebesgue measure on $\Reals^{|\T_*|}$.

Let us now recall the definition of the clusters $(H_\alpha)_{\alpha\in\A}$ and their Gibbs weights $(V_\alpha)_{\alpha\in\A}$ in the paragraph above equation (\ref{Vs}). Suppose that the cardinality of $\LL(\T)$ is equal to $m$. Let us look at all possible choices of $m$ pure states $H_{\alpha_t}$ for $t\in \LL(\T)$ indexed by the leaves $\alpha_t\in \LL(\A)=\Natural^r$ that ``form the same pattern" according to their overlaps as the tree $\T$. More precisely, we will denote $\bar{\alpha} := (\alpha_t)_{t\in \LL(\T)}$ and consider the set  
$$
\A(\CC) =
\Bigl\{ \bar{\alpha} \in \LL(\A)^m \ \bigr| \
\alpha_t \wedge \alpha_{t'} = t\wedge t' \ \mbox{ for all }\ t,t'\in \LL(\T)
\Bigr\}.
$$
Then it should be obvious that the event $\OO = \OC$ defined in (\ref{OCset}) can be written as a disjoint union $\OO = \bigcup_{\bar{\alpha}\in \A(\CC)} \OO(\bar{\alpha})$, where (recall the definition of $\RR(t)$ in (\ref{Rt}))
$$
\OO(\bar{\alpha}) = \Bigl\{ (\sigma^1,\ldots,\sigma^n) \ \bigr|\
\sigma^\ell\in H_{\alpha_t} \mbox{ for all } t\in\LL(\T) \mbox{ and } \ell \in \RR(t) \Bigr\}.
$$
Then, we can write
$$
\p_\OO(W\in B) = \frac{\e\bigl\la \Ind(W\in B) \Ind_{\OO}\bigr\ra}{\p(\OO)}
= \frac{1}{\p(\OO)}\e \sum_{ \bar{\alpha} \in \A(\CC)} \bigl\la \Ind(W\in B) \Ind_{\OO(\bar{\alpha})}\bigr\ra. 
$$
On the event $\OO(\bar{\alpha})$, the vector of weights $W=(W_t)_{t\in\T_*}$  can also be written as a vector of cluster weights $V_\alpha$ in (\ref{Vs}) indexed by the vertices $\alpha$ in the subtree formed by all paths from the root to the leaves $(\alpha_t)_{t\in \LL(\T)}$. Let us call this vector $V(\bar{\alpha})$. Also, obviously,
$$
\bigl\la \Ind_{\OO(\bar{\alpha})}\bigr\ra = \prod_{t\in\LL(\T)} V_{\alpha_t}^{|\RR(t)|}
$$
and, therefore,
$$
\p_\OO(W\in B) =
\frac{1}{\p(\OO)}\e \sum_{\bar{\alpha} \in \A(\CC)} \Ind\bigl(V(\bar{\alpha}) \in B\bigr) 
\prod_{t\in\LL(\T)} V_{\alpha_t}^{|\RR(t)|}.
$$
To finish the proof of Lemma \ref{LemAC}, it is enough to show that the distribution of $V(\bar{\alpha})$ is absolutely continuous with respect to the Lebesgue measure. For the remainder of this section, we will forget about the configuration $\CC$ and will focus on proving the absolute continuity of the distribution of cluster weights $(V_\alpha)_{\alpha\in F}$ indexed by an arbitrary finite subset $F$ of $\A\setminus\{*\}$. Of course, this will be based on the properties of the Ruelle probability cascades (RPC), so we will first recall the construction of these cascades and how it relates to the weights $V_\alpha$.

Recall the sequence of parameters in (\ref{zetas}). For each $\alpha\in \A\setminus \Natural^r$, let $\Pi_\alpha$ be a Poisson process on $(0,\infty)$ with the mean measure $\zeta_{p}x^{-1-\zeta_{p}}dx$ with $p=|\alpha|$, and we assume that these processes are independent for all $\alpha$. Let us arrange all the points in $\Pi_\alpha$ in the decreasing order,
\begin{equation}
u_{\alpha 1} > u_{ \alpha 2} >\ldots >u_{\alpha n} > \ldots,
\label{us}
\end{equation}
and enumerate them using the children $(\alpha n)_{n\geq 1}$ of the vertex $\alpha$. Given a vertex $\alpha\in \A\setminus \{*\}$ and the path $p(\alpha)$ in (\ref{pathtoleaf}), we define
\begin{equation}
w_\alpha = \prod_{\beta \in p(\alpha)} u_{\beta},
\label{ws}
\end{equation}
and for the leaf vertices $\alpha \in \LL(\A) = \Natural^r$ we define
\begin{equation}
v_\alpha = \frac{w_\alpha}{\sum_{\beta\in \Natural^r} w_\beta}.
\label{vs}
\end{equation}
For other vertices $\alpha\in \A\setminus \LL(\A)$ we define
\begin{equation}
v_\alpha = \sum_{\beta\in \LL(\A),\,\beta\succ \alpha} v_\beta.
\label{vsall}
\end{equation}
Of course, this definition implies that $v_\alpha = \sum_{n\geq 1} v_{\alpha n}$ when $|\alpha|<r$.
Notice that, for a given $\alpha$, the sequence of weights $(v_{\alpha n})_{n\geq 1}$ is not necessarily decreasing. For example, when $r=2$, sequences $(u_n)_{n\geq 1}$ and $(u_{nm})_{m\geq 1}$ for all $n$ are decreasing by construction, but $v_n$ is proportional to $u_n\sum_{m\geq 1} u_{nm}$ and does not have to be decreasing. Let us now rearrange the vertex labels so that the weights indexed by children will be decreasing. For each $\alpha\in \A\setminus \Natural^r$, let $\pi_\alpha: \Natural \to \Natural$ be a bijection such that the sequence $(v_{\alpha \pi_\alpha(n)})_{n\geq 1}$ is decreasing. Using these ``local rearrangements" we define a global bijection $\pi: \A\to \A$ in a natural way, as follows. We let $\pi(*)=*$ and then define
\begin{equation}
\pi(\alpha n) = \pi(\alpha) \pi_{\pi(\alpha)}(n)
\label{permute}
\end{equation}
recursively from the root to the leaves of the tree. Finally, we define
\begin{equation}
V_\alpha = v_{\pi(\alpha)} \ \mbox{ for all }\  \alpha\in \A.
\label{Vs2}
\end{equation}
It is not a coincidence that we used here the same notation as in (\ref{Vs}), since they have the same distribution. This relationship between cluster weights of a random measure $G$ and the RPC is a well-known consequence of the Ghirlanda-Guerra identities (see Section 2.4 in \cite{SKmodel}). Therefore, our goal is to prove the following.

\begin{lemma}\label{Lem3}
The distribution of weights $(V_\alpha)_{\alpha\in F}$ in (\ref{Vs2}) indexed by an arbitrary finite subset $F$ of $\A\setminus\{*\}$ is absolutely continuous with respect to the Lebesgue measure on $\Reals^{|F|}$.
\end{lemma}
Let us first introduce some more notation and recall some definitions. Let $(u_n)_{n\geq 1}$ be the decreasing enumeration of a Poisson process on $(0,\infty)$ with the mean measure $x u^{-1-x}du$ for some $x\in (0,1)$ and let 
\begin{equation}
U = \sum_{n\geq 1} u_n \ \mbox{ and }\
p_n = \frac{u_n}{U} \ \mbox{ for }\ n\geq 1.
\end{equation}
The distribution of the sequence $(p_n)_{n\geq 1}$ is called the \emph{Poisson-Dirichlet distribution} $PD(x)$ (or $PD(x,0)$). It is well known that the distribution of finitely many coordinates of $(p_n)$ is absolutely continuous. For example, Proposition 47 in \cite{Pitman} gives some representation for the density, but the existence of the density is also easy to see directly from the representation of this process in Proposition 8 in \cite{Pitman}.

Let us consider $a<x$. Then the distribution of $(p_n)_{n\geq 1}$ under the change of density 
$U^a/\e U^a$ is called the \emph{Poisson-Dirichlet distribution} $PD(x,-a)$. The usual condition $a<x$ ensures that $\e U^a <\infty$ and the change of density is well defined (see e.g. Lemma 2.1 in \cite{SKmodel}). The definition of this distribution in Section 1.1 in \cite{Pitman} was different but its equivalence to this one was shown in Proposition 14 there. (In \cite{Pitman}, the parameter $-a$ was denoted $\theta$ and the condition was stated as $\theta>-x$.) It is easy to see that the distribution of finitely many coordinates of $(p_n)$ under $PD(x,-a)$ is also absolutely continuous. Indeed, for any $N\geq 1$ and a measurable set $A$ in $\Reals^N$ of Lebesgue measure $0$, by H\"older's inequality,
\begin{equation}
\e U^a \Ind\bigl((p_n)_{n\leq N}\in A\bigr)
\leq
(\e U^{a(1+\eps)})^{1/(1+\eps)} \p\bigl((p_n)_{n\leq N}\in A\bigr)^{\eps/(1+\eps)}=0,
\label{absxa}
\end{equation}
for small enough $\eps>0$ such that $a(1+\eps)<x$, in which case $\e U^{a(1+\eps)}<\infty$.

For each $\alpha\in \Natural^{r-1}$, let us now consider the sequence
\begin{equation}
p_{\alpha n} = \frac{V_{\alpha n}}{V_\alpha} \ \mbox{ for }\ n\geq 1.
\label{pa}
\end{equation}
By definition, this sequence is decreasing and $\sum_{n\geq 1} p_{\alpha n}=1$.  The following holds.
\begin{lemma}\label{Lem4}
For each $\alpha\in \Natural^{r-1}$, the sequence $(p_{\alpha n})_{n\geq 1}$ in (\ref{pa}) has distribution $PD(\zeta_{r-1},-\zeta_{r-2})$. These sequences are independent of each other and of $(V_{\alpha})_{|\alpha|\leq r-1}$.
\end{lemma}
First, let us show how this implies Lemma \ref{Lem3}.

\medskip
\noindent
\textbf{Proof of Lemma \ref{Lem3}.} This now follows easily by induction on $r$. For $r=1$, this is just absolute continuity of weights from the Poisson-Dirichlet distribution $PD(\zeta_0)$. To make an induction step, we use a well-known fact that the array $(V_\alpha)_{|\alpha| \leq r-1}$ can be constructed as in (\ref{us}) -- (\ref{Vs2}) with $r$ replaced by $r-1$ and $\zeta_{r-1}$ removed from the sequence (\ref{zetas}). This observation goes back to \cite{Ruelle}, but is also a trivial consequence of the Ghirlanda-Guerra identities. (In any case, the proof of this fact will appear below as a byproduct of the proof of Lemma \ref{Lem4}.) By the induction hypothesis, this implies that the distribution of finitely many nonroot coordinates of $(V_\alpha)_{|\alpha| \leq r-1}$ is absolutely continuous. To include coordinates $V_{\alpha n}$ for $\alpha\in\Natural^{r-1}$ and $n\geq 1$, we write them as $V_{\alpha n} = V_\alpha p_{\alpha n}$ and use Lemma \ref{Lem4} together with the observation in (\ref{absxa}) about absolutely continuity of the distribution of finitely many coordinates under $PD(x,-a)$.
\qed

\medskip
\noindent
\textbf{Proof of Lemma \ref{Lem4}.} 
We only need to consider the case $r\geq 2$. For each $\alpha\in\Natural^{r-2}$, consider the process $(u_{\alpha n}, (u_{\alpha n m})_{m\geq 1} )_{n\geq 1}$ and let $U_{\alpha n} := \sum_{m\geq 1} u_{\alpha n m}$. If we define
$$
d_{\alpha n m} = \frac{v_{\alpha n m}}{v_{\alpha n}} =  \frac{u_{\alpha n m}}{U_{\alpha n}}
$$
then $Y_{\alpha n} := (d_{\alpha n m})_{m\geq 1}$ has the Poisson-Dirichlet distribution $PD(\zeta_{r-1})$. Notice that the random variables $(U_{\alpha n}, Y_{\alpha n})_{n\geq 1}$ are i.i.d. and independent of $(u_{\alpha n})_{n\geq 1}$. Moreover, all these processes are independent over $\alpha\in \Natural^{r-2}$, and also independent of $U_{r-2} = (u_\alpha)_{|\alpha|\leq r-2}$. 

For a fixed $\alpha\in\Natural^{r-2}$, let $\pi_\alpha:\Natural \to \Natural$ be a bijection such that the sequence $(u_{\alpha \pi_\alpha(n)} U_{\alpha \pi_\alpha(n)})_{n\geq 1}$ is decreasing. This is exactly the same permutation defined in the paragraph above (\ref{permute}) since, for a fixed $\alpha\in\Natural^{r-2}$, $v_{\alpha n}$ is proportional to $u_{\alpha n}U_{\alpha n}$. Since $(u_{\alpha n})_{n\geq 1}$ is a Poisson process with the mean measure $\zeta_{r-2} \,x^{-1-\zeta_{r-2}} dx$, Theorem 2.6 in \cite{SKmodel} (Proposition A.2 in \cite{Bolthausen}) implies that 
\begin{equation}
\bigl(u_{\alpha \pi_\alpha(n)} U_{\alpha \pi_\alpha(n)} , Y_{\alpha \pi_\alpha(n)}\bigr)_{n\geq 1}
\stackrel{d}{=}
\bigl(u_{\alpha n} c, Y_{\alpha n}'\bigr)_{n\geq 1}
\label{Lem4inproof1}
\end{equation}
where $c= \bigl(\e U_{\alpha 1}^{\zeta_{r-2}} \bigr)^{1/\zeta_{r-2}}$, $(u_{\alpha n})_{n\geq 1}$ and $(Y_{\alpha n}')_{n\geq 1}$ on the right hand side are independent, and the random variables $(Y_{\alpha n}')_{n\geq 1}$ are i.i.d. with the distribution of $Y_{\alpha 1} = (d_{\alpha 1 m})_{m\geq 1}$ under the change of density  
$$
U_{\alpha 1}^{\zeta_{r-2}} \bigr/ \,\e U_{\alpha 1}^{\zeta_{r-2}},
$$
which is precisely the Poisson-Dirichlet distribution $PD(\zeta_{r-1},-\zeta_{r-2})$. It remains to notice that the weights $(V_{\alpha})_{|\alpha|\leq r-1}$ are, obviously, a function of the arrays 
\begin{equation}
U_{r-2} = (u_\alpha)_{|\alpha|\leq r-2} \ \mbox{ and }\ 
\bigl(u_{\alpha \pi_\alpha(n)} U_{\alpha \pi_\alpha(n)} \bigr)_{\alpha\in \Natural^{r-2},n\geq 1}
\label{Lem4inproof2}
\end{equation}
and are, therefore, independent of the random variables $Y_{\alpha \pi_\alpha(n)}$, which are i.i.d. for all $\alpha\in\Natural^{r-2}$ and $n\geq 1$ and have the distribution $PD(\zeta_{r-1},-\zeta_{r-2})$. In particular, the permutation $\pi$ defined in (\ref{permute}), restricted to $|\alpha|\leq r-1$, will be a function of these arrays and, therefore, 
$$
\bigl(d_{\pi(\alpha n)m} \bigr)_{m\geq 1} =
Y_{\pi(\alpha n)}  = Y_{\pi(\alpha) \pi_{\pi(\alpha)}(n)} 
$$
are still i.i.d. over all $\alpha\in\Natural^{r-2}$ and $n\geq 1$, have distribution $PD(\zeta_{r-1},-\zeta_{r-2})$, and are independent of $(V_{\alpha})_{|\alpha|\leq r-1}$. This finishes the proof since, by the definition (\ref{pa}), for $\alpha n\in \Natural^{r-1}$, 
$$
p_{\alpha n m} =
 \frac{V_{\alpha n m}}{V_{\alpha n}}
=  \frac{v_{\pi(\alpha n) m}}{v_{\pi(\alpha n)}}
= d_{\pi(\alpha n)m}.
$$
Finally, let us notice that the above argument also proves the fact mentioned in the proof of Lemma \ref{Lem3}, namely, that the array $(V_\alpha)_{|\alpha| \leq r-1}$ can be constructed as in (\ref{us}) -- (\ref{Vs2}) with $r$ replaced by $r-1$ and $\zeta_{r-1}$ removed from the sequence (\ref{zetas}). This is because $(V_\alpha)_{|\alpha| \leq r-1}$ is constructed from the arrays in (\ref{Lem4inproof2}) as in (\ref{us}) -- (\ref{Vs2}) and, by (\ref{Lem4inproof1}), for each $\alpha\in \Natural^{r-2}$, the second array in (\ref{Lem4inproof2}) is, up to a  factor $c$, a Poisson process with the mean measure $\zeta_{r-2} \,x^{-1-\zeta_{r-2}} dx$. Of course, this constant factor $c$ will cancel at the step (\ref{vs}), so the claim follows.
\qed

\section{From replicas to the Gibbs measure} \label{Sec5label}

In this section, we will show how Theorem \ref{Th1} can be deduced from Theorem \ref{Th3}. The main idea is that when the sample size $n$ goes to infinity, there will be many replicas in any given subset of pure states, and the statement in Theorem \ref{Th3} about spins and cluster weights corresponding to the sample can be translated into a statement in Theorem \ref{Th1} about spins inside pure states and cluster weights of the Gibbs measure.

Before we begin the proof, let us first notice that Theorem \ref{Th1} follows from its analogue for finite subsets of the tree $\A$, as follows. Let us consider integers $d\geq 1$ and $N\geq 1$ that will be fixed throughout this section (note that now the notation $N$ is not related to the number of coordinates of the system in the introduction). Let $[d]=\{1,\ldots, d\}$ and let
$$
\A_d = \{*\}\cup [d] \cup [d]^2 \cup\ldots\cup [d]^r \subseteq \A
$$
be a $d$-regular subtree of $\A$. Any finite subset of $\A$ will be covered by $\A_d$ for $d$ sufficiently large (depending on the subset). Now, recall the array $S_\alpha = (S(\sigma^{\alpha n}))_{n\geq 1}$ in (\ref{Salpha}) and let us truncate it to the array
\begin{equation}
S_{\alpha,N} = \bigl(S(\sigma^{\alpha n}) \bigr)_{n\leq N}
\label{SalphaDN}
\end{equation}
generated by a sample $(\sigma^{\alpha n})_{n\leq N}$ of size $N$ from the pure state $H_{\alpha}$. We will only consider these arrays for $\alpha \in [d]^r = \LL(\A_d)$, so we will need to restrict the notion of hierarchical exchangeability to the finite tree $\A_d$. Similarly to (\ref{setH}), let
\begin{equation}
{\cal H}_d = \Bigl\{ \pi : [d]^r\to [d]^r \ \bigr|\ \pi \mbox{ is a bijection}, \pi(\alpha)\wedge \pi(\beta) = \alpha\wedge\beta \mbox{ for all } \alpha,\beta\in [d]^r \Bigr\}.
\label{setHDN}
\end{equation}
Then, naturally, we will call a finite array $(X_\alpha)_{\alpha\in [d]^r}$ {hierarchically exchangeable} if 
\begin{equation}
\bigl(X_{\pi(\alpha)} \bigr)_{\alpha\in [d]^r}
\stackrel{d}{=}
\bigl(X_\alpha \bigr)_{\alpha\in [d]^r}
\label{HexchDFDN}
\end{equation}
for all $\pi\in {\cal H}_d$. It is obvious that, in order to prove Theorem \ref{Th1}, it is sufficient to show the following for all $d,N \geq 1$.
\begin{thmref}{Th1}
The array of spins $(S_{\alpha,N})_{\alpha\in [d]^r}$ defined in (\ref{SalphaDN}) is hierarchically exchangeable and independent of the array of cluster weights $(V_\alpha)_{\alpha\in\A_d\setminus\{*\}}$.
\end{thmref}
To prove this, we will apply Theorem \ref{Th3} to the following set of configurations $\CC=(\T,\PP)$,
\begin{equation}
\CC(n,d,N) = \Bigl\{ \CC=(\T,\PP) \ \bigr|\   \A_d\subseteq \T \mbox{ and } |\PP^{-1}(t)|\geq N \mbox{ for } t\in [d]^r \Bigr\}
\label{CndN}
\end{equation}
(this set depends on $n$ through the mapping $\T$). In words, the tree $\T$ contains $\A_d$ (so it is big enough) and at least $N$ replica indices are mapped by $\PP$ into each leaf $t\in [d]^r = \LL(\A_d)\subseteq \LL(\T)$. For a given configuration $\CC\in \CC(n,d,N)$ and $t\in [d]^r$, let $\RR_N(t)$ be the set of the smallest $N$ replica indices in $\PP^{-1}(t)$ (we choose the smallest $N$ just for certainty, and arbitrary $N$ would do) and define $\RR_{d,N} = \bigcup_{t\in [d]^r} \RR_N(t)$. Let us recall the definition of $S^n$ in (\ref{SpinsEll}) and (\ref{Sn}) and, similarly, define
\begin{equation}
S^{d,N} = \bigl(S(\sigma^\ell) \bigr)_{\ell\in \RR_{d,N}}.
\label{SndN}
\end{equation}
In other words, we are now only interested in a set of $N$ replicas for each of the leaves in $[d]^r.$ Similarly to (\ref{OCset}), let us define the event
\begin{equation}
\OO(\CC,d,N) = \Bigl\{ (\sigma^\ell)_{\ell\in \RR_{d,N}} \ \bigr|\ \sigma^\ell\cdot\sigma^{\ell'}\in I_{\PP(\ell)\wedge \PP(\ell')}  \mbox{ for all } \ell,\ell' \in \RR_{d,N}\Bigr\},
\label{OCdN}
\end{equation}
which involves only the replicas with indices in $\RR_{d,N}$ and, similarly to the definition of $\p_{\OC}$ in (\ref{PO}), we let
\begin{equation}
\p_{\OO(\CC,d,N)}(\ \cdot \ ) = \frac{\p(\ \cdot\ \cap \OO(\CC,d,N))}{\p(\OO(\CC,d,N))}.
\label{POdN}
\end{equation}
We will need the following simple consequence of the Ghirlanda-Guerra identities (\ref{GG}). 
\begin{lemma}\label{Lem5}
For any $\CC\in \CC(n,d,N)$, we have 
\begin{equation}
\p_{\OC}\bigl(S^{d,N}\in \ \cdot \ \bigr) = \p_{\OO(\CC,d,N)}\bigl(S^{d,N}\in \ \cdot \ \bigr).
\label{Lem5eq}
\end{equation}
\end{lemma}
\textbf{Proof.} Consider any replica index in $\{1,\ldots,n\}\setminus\RR_{d,N}$ not appearing in $S^{d,N}$. For simplicity of notation, suppose that this index is $n$. Let
$$
p=\max_{\ell<n}\PP(n)\wedge\PP(\ell)
$$
and let
$$
\OC^- = \Bigl\{ (\sigma^\ell)_{1\leq\ell\leq n-1}\ \bigr|\
\sigma^\ell\cdot\sigma^{\ell'}\in I_{\PP(\ell)\wedge\PP(\ell')}
\mbox{ for all }1\leq\ell,\ell'\leq n-1\Bigr\}.
$$
We will show that, for any bounded measurable function $f$ of the retained spins $S^{d,N}$,
$$
\e\la f\Ind_{\OC}\ra=c\,\e\la f\Ind_{\OC^-}\ra,
$$
where $c>0$ depends only on the prescribed configuration.

Suppose first that $p<r$. Let
$$
J=\{\ell<n:\PP(\ell)\wedge\PP(n)\geq p\},\qquad M=|J|.
$$
Partition $J$ into the nonempty sets $J_1,\ldots,J_b$ corresponding to the distinct depth-$(p+1)$ children occupied by these replicas. Write $M_j=|J_j|$ and choose $i_j\in J_j$. Replica $n$ must belong to their common depth-$p$ cluster, but to none of these $b$ occupied children. By ultrametricity, almost surely,
$$
\Ind_{\OC}=\Ind_{\OC^-}\Bigl(
\Ind(d_I(R_{i_1,n})\geq p)
-\sum_{j=1}^b\Ind(d_I(R_{i_j,n})\geq p+1)\Bigr).
$$
Indeed, on $\OC^-$ the events in the sum are disjoint subsets of the first event, and the overlap constraints with replicas outside $J$ follow automatically. When $p=0$, the first indicator is identically one.

Since $\zeta\{x:d_I(x)\geq p\}=1-\zeta_{p-1}$, the Ghirlanda-Guerra identities give
$$
\e\la f\Ind_{\OC^-}\Ind(d_I(R_{i_1,n})\geq p)\ra
=\frac{M-\zeta_{p-1}}{n-1}\e\la f\Ind_{\OC^-}\ra.
$$
Similarly, for each $j\leq b$,
$$
\e\la f\Ind_{\OC^-}\Ind(d_I(R_{i_j,n})\geq p+1)\ra
=\frac{M_j-\zeta_p}{n-1}\e\la f\Ind_{\OC^-}\ra.
$$
Subtracting and using $\sum_{j=1}^bM_j=M$, we obtain
$$
\e\la f\Ind_{\OC}\ra
=\frac{b\zeta_p-\zeta_{p-1}}{n-1}\e\la f\Ind_{\OC^-}\ra.
$$
The coefficient is strictly positive by (\ref{zetas}).

If $p=r$, replica $n$ belongs to an already occupied terminal cluster. Let $M$ be the number of remaining replicas assigned to that leaf and choose one of their indices $i$. In this case, ultrametricity gives
$$
\Ind_{\OC}=\Ind_{\OC^-}\Ind(d_I(R_{i,n})\geq r),
$$
and the Ghirlanda-Guerra identities imply
$$
\e\la f\Ind_{\OC}\ra
=\frac{M-\zeta_{r-1}}{n-1}\e\la f\Ind_{\OC^-}\ra.
$$
This coefficient is also strictly positive.

For completeness, when $n=2$, the one-replica identity used in these calculations follows from (\ref{GG}) with two replicas. If $f$ depends only on $S^1$, replica exchangeability gives
$$
\e\la f\psi(R_{1,2})\ra
=\e\la f\psi(R_{1,3})\ra
=\frac12\e\la f\ra\,\e\la\psi(R_{1,2})\ra
+\frac12\e\la f\psi(R_{1,2})\ra,
$$
and hence $\e\la f\psi(R_{1,2})\ra=\e\la f\ra\,\e\la\psi(R_{1,2})\ra$.

Applying the calculation with $f=\Ind(S^{d,N}\in A)$ and with $f=1$, the same positive coefficient cancels between numerator and denominator. Therefore,
$$
\p_{\OC}(S^{d,N}\in A)=\p_{\OC^-}(S^{d,N}\in A).
$$
We can now remove the other replica indices outside $\RR_{d,N}$ one by one in the same way. This proves (\ref{Lem5eq}).
\qed

\medskip
\medskip
\noindent
\textbf{Remark.} Notice that the right hand side of (\ref{Lem5eq}) does not really depend on the configuration $\CC$ since the set $\RR_{d,N}$ involves $N$ replicas assigned to the leaves $[d]^r$ of the tree $\A_d$, and we can relabel those replicas using indices $1,\ldots, N d^r.$ Let $\CC_{d,N}$ be a configuration consisting of the tree $\A_d$ and a map $\PP_{d,N}$ that maps exactly $N$ indices in $\{1,\ldots, N d^r\}$ to each leaf in $[d]^r$. Then the equation (\ref{Lem5eq}) can be rewritten as
\begin{equation}
\p_{\OC}\bigl(S^{d,N}\in \ \cdot \ \bigr) = \p_{\OO(\CC_{d,N})  }\bigl(S^{d,N}\in \ \cdot \ \bigr).
\label{Lem5eq2}
\end{equation} 
We use the same notation $S^{d,N}$ on the right hand side but, of course, we need to change the definition of $S^{d,N}$ to take into account this relabeling of indices. In fact, for clarity, let us index the $N$ replicas mapped by $\PP_{d,N}$ into the leaf $\alpha\in [d]^r = \LL(\A_d)$ by $\sigma^{(\alpha, 1)},\ldots, \sigma^{(\alpha, N)}.$ Then $S^{d,N}$ on the right hand side of (\ref{Lem5eq2}) is understood as 
\begin{equation}
S^{d,N} = \bigl(S(\sigma^{(\alpha,\ell)}) \bigr)_{\alpha\in [d]^r,\ell \leq N}.
\label{SndNred}
\end{equation}
Notice that we use the notation $\sigma^{(\alpha,\ell)}$ here to distinguish these (usual, unconditional) replicas from the Gibbs measure $G$ from the replicas $\sigma^{\alpha \ell}$ in (\ref{Salpha}), which denoted the sample from conditional Gibbs measure $G_\alpha$ on the pure state $H_\alpha$.
\qed

\medskip
\medskip
\noindent
For a given configuration $\CC=(\T,\PP)$, let us recall the definition of $W = (W_t)_{t\in \T_*}$ in (\ref{WC}) and (\ref{WCT}), which represent the cluster weights around the sample on the event $\OC$. For a configuration $\CC \in \CC(n,d,N)$ in (\ref{CndN}), we will denote by 
\begin{equation}
W^d = (W_t)_{t\in \A_d\setminus\{*\}}
\label{WndN}
\end{equation}
the subset of these weights along the subtree $\A_d\subseteq \T$. Let us recall the definition of the sample configuration $\CC_n = (\T_n,\PP_n)$ in (\ref{sampleconfig}) and consider two events
\begin{align}
\EE_1(n)
& =
\bigcup_{\CC\in  \CC(n,d,N)} \Bigl\{S^{d,N} \in A, W^d\in B, \CC_n = \CC \Bigr\}, 
\label{E1}
\\
\EE_2(n)
&=
\bigcup_{\CC\in  \CC(n,d,N)} \Bigl\{W^d\in B, \CC_n = \CC \Bigr\}.
\label{E2}
\end{align}
To understand what these events represent, let us see what they will look like with high probability when the sample size $n\to\infty$. When $n$ gets large, with high probability, at least $N$ replicas will fall into each of the pure states $H_\alpha$ for $\alpha\in [d]^r$. First of all, this means that with high probability the sample configuration $\CC_n \in  \CC(n,d,N)$. Second, conditionally on this event that at least $N$ replicas fall into each of the pure states $H_\alpha$ for $\alpha\in [d]^r$, what are $S^{d,N}$ and $W^d$ in (\ref{E1}) and (\ref{E2})?
Recall that $\CC_n = \CC$ means that the event $\WC$ in (\ref{WCset}) occurs and, for each vertex $t\in \T\setminus \LL(\T)$, the cluster weights indexed by its children are arranged in decreasing order. The pure states $H_\alpha$ and the weights $V = (V_\alpha)_{\alpha\in \A}$ in (\ref{Vs}) of the clusters around the pure states were labelled in a similar fashion in (\ref{purelabels}). This implies that whenever $\CC_n = \CC \in \CC(n,d,N)$ and at least $N$ replicas fall into each of the pure states $H_\alpha$ for $\alpha\in [d]^r$, we must have $W^d = (V_\alpha)_{\alpha\in\A_d\setminus\{*\}}.$ Moreover, in this case, the spins $S^{d,N}$ correspond to $N$ replicas sampled from each of the pure states $H_\alpha$ for $\alpha\in [d]^r$, i.e. $S^{d,N} = (S_{\alpha,N})_{\alpha\in [d]^r}$ defined in (\ref{SalphaDN}). This implies that
\begin{align}
\lim_{n\to\infty}\p\bigl(\EE_1(n)\bigr)
& =
\p\bigl((S_{\alpha,N})_{\alpha\in [d]^r} \in A, (V_\alpha)_{\alpha\in\A_d\setminus\{*\}} \in B\bigr),
\label{E1lim}
\\
\lim_{n\to\infty}\p\bigl(\EE_2(n)\bigr)
& =
\p\bigl((V_\alpha)_{\alpha\in\A_d\setminus\{*\}} \in B\bigr).
\label{E2lim}
\end{align}
To finish the proof of Theorem \ref{Th1}${}^\prime$, it remains to show the following.
\begin{lemma} \label{Lem6}
We have, 
\begin{equation}
\p\bigl(\EE_1(n)\bigr) = \p_{\OO(\CC_{d,N})  }\bigl(S^{d,N}\in A \bigr) \p\bigl(\EE_2(n)\bigr).
\end{equation}
\end{lemma}

\smallskip
\noindent
\textbf{Proof.} First of all, when we defined the sample configuration $\CC_n$ in (\ref{sampleconfig}) we explained that the events $\CC_n =\CC$ are disjoint for different $\CC$ and $\{\CC_n =\CC\} = \WC \cap \OC$. Therefore,
\begin{align*}
\p\bigl(\EE_1(n) \bigr)
&=
\sum_{\CC\in  \CC(n,d,N)}  \p\bigl( \{S^{d,N} \in A\} \cap \{W^d\in B\}\cap \WC\cap \OC\bigr),
\\
\p\bigl(\EE_2(n) \bigr)
&=
\sum_{\CC\in  \CC(n,d,N)}  \p\bigl(\{W^d\in B\}\cap \WC\cap \OC\bigr).
\end{align*}
Notice that $\{W^d\in B\}\cap \WC$ is an event which involves only the weights $W = (W_t)_{t\in \T_*}$ and can be written as $\{W\in B'\}$ for some set $B'$. Therefore, Theorem \ref{Th3} implies that 
\begin{align*}
&
\p\bigl( \{S^{d,N} \in A\} \cap \{W^d\in B\}\cap \WC\cap \OC\bigr)
\\
&\hspace{2.2cm}= 
\p_{\OC}\bigl(S^{d,N} \in A)\, \p\bigl(\{W^d\in B\}\cap \WC\cap \OC\bigr).
\end{align*}
Finally, using (\ref{Lem5eq2}), we can write
\begin{align*}
\p\bigl(\EE_1(n) \bigr)
&=
\p_{\OO(\CC_{d,N})  }\bigl(S^{d,N}\in A \bigr)
\sum_{\CC\in  \CC(n,d,N)}  \p\bigl(\{W^d\in B\}\cap \WC\cap \OC\bigr)
\\
&=
\p_{\OO(\CC_{d,N})  }\bigl(S^{d,N}\in A \bigr)\,  \p\bigl(\EE_2(n) \bigr),
\end{align*}
which finishes the proof.
\qed
 
\medskip
\medskip
\noindent
Together with (\ref{E1lim}) and (\ref{E2lim}), Lemma \ref{Lem6} implies 
$$
\p\bigl((S_{\alpha,N})_{\alpha\in [d]^r} \in A, (V_\alpha)_{\alpha\in\A_d\setminus\{*\}} \in B\bigr)
=
\p_{\OO(\CC_{d,N})  }\bigl(S^{d,N}\in A \bigr)\,
\p\bigl((V_\alpha)_{\alpha\in\A_d\setminus\{*\}} \in B\bigr).
$$
Therefore, $(S_{\alpha,N})_{\alpha\in [d]^r}$ and $(V_\alpha)_{\alpha\in\A_d\setminus\{*\}}$ are independent  and, recalling (\ref{SndNred}),
$$
\p\bigl((S_{\alpha,N})_{\alpha\in [d]^r} \in A \bigr)
=
\p_{\OO(\CC_{d,N})  }\bigl(\bigl(S(\sigma^{(\alpha,\ell)}) \bigr)_{\alpha\in [d]^r,\ell \leq N} \in A \bigr).
$$
The hierarchical exchangeability of $(S_{\alpha,N})_{\alpha\in [d]^r}$ follows, because of the obvious invariance of the event $\OO(\CC_{d,N})$ under the permutations $\pi\in {\cal H}_d$ in (\ref{setHDN}),
$$
\bigl(\sigma^{(\alpha,\ell)}\bigr)_{\alpha\in [d]^r,\ell \leq N} \in \OO(\CC_{d,N})
\Longleftrightarrow
\bigl(\sigma^{(\pi(\alpha),\ell)}\bigr)_{\alpha\in [d]^r,\ell \leq N} \in \OO(\CC_{d,N}).
$$
This finishes the proof of Theorem \ref{Th1}${}^\prime$ and, thus, Theorem \ref{Th1}.
\qed

\end{document}